\documentclass[11pt]{article}
\usepackage[dvipsnames,table,xcdraw]{xcolor}
\usepackage{amsmath}
\usepackage{amsthm}
\usepackage{mathtools}
\usepackage{amssymb}
\usepackage{graphicx}
\usepackage[hidelinks]{hyperref}
\usepackage{float}
\usepackage{mathtools}
\usepackage{makecell}
\usepackage{hhline}
\usepackage{nicematrix}
\usepackage{array}
\usepackage{subcaption}

\usepackage{tkz-berge}
\usetikzlibrary{arrows.meta}
\def\niceArrow{-{Stealth[length=2.5mm]}}

\usepackage{tikz}
\usetikzlibrary{positioning,decorations.pathreplacing,arrows.meta,calc}

\usetikzlibrary{matrix, fit, backgrounds, positioning}

\usepackage[margin=16pt]{caption}

\usepackage[a4paper, total={6.25in, 9.49in}]{geometry}

\usepackage[textsize=tiny]{todonotes}

\def\R{\mathbb R}
\def\Z{\mathbb Z}

\numberwithin{equation}{section}

\newcommand{\nat}{^\natural}

\newcommand{\FBOX}{\hspace*{\fill}$\rule{0.17cm}{0.17cm}$}

\newcommand{\ol}{\overline}

\newtheorem{THM}{THEOREM}[section]

\newtheorem{theorem}[THM]{Theorem}
\newtheorem{lemma}[THM]{Lemma}
\newtheorem{claim}[THM]{Claim}
\newtheorem{corollary}[THM]{Corollary}
\newtheorem{proposition}[THM]{Proposition}

\theoremstyle{definition}
\newtheorem{remark}{Remark}[section]

\usepackage{titling}
\usepackage[hang,flushmargin]{footmisc}
\usepackage[yyyymmdd]{datetime}

\newcommand\nnfootnote[1]{%
  \begin{NoHyper}
    \renewcommand\thefootnote{}\footnote{#1}%
    \addtocounter{footnote}{-1}%
  \end{NoHyper}
}

\begin{document}

\title{Prefix-bounded matrices
}

\author{N\'ora A.\ Borsik \thanks{Department of Operations Research, E\"otv\"os Lor\'and University Budapest, P\'azm\'any P.\ s.\ 1/c, Budapest H-1117. E-mail: {\tt nborsik@gmail.com}}
  \and Andr\'as Frank \thanks{HUN-REN--ELTE Egerv\'ary Research Group on Combinatorial Optimization, Department of Operations Research, E\"otv\"os Lor\'and University Budapest, P\'azm\'any P.\ s.\ 1/c, Budapest H-1117. E-mail: {\tt andras.frank@ttk.elte.hu}}
  \and P\'eter Madarasi \thanks{HUN-REN Alfr\'ed R\'enyi Institute of Mathematics, and Department of Operations Research, E\"otv\"os Lor\'and University Budapest, P\'azm\'any P.\ s.\ 1/c, Budapest H-1117. E-mail: {\tt madarasip@staff.elte.hu}}
  \and Tam\'as Tak\'acs \thanks{Department of Operations Research, E\"otv\"os Lor\'and University Budapest, P\'azm\'any P.\ s.\ 1/c, Budapest H-1117. E-mail: {\tt takacst2003@gmail.com}}
}

\date{\vspace*{-24pt}}

\maketitle
\nnfootnote{Corresponding author: P\'eter Madarasi}

\begin{abstract}
  By unifying various earlier extensions of alternating sign matrices (ASMs), we introduce the notion of prefix-bounded matrices (PBMs).
  It is shown that the convex hull of these matrices forms the intersection of two special generalized polymatroids.
  This implies --- in a more general form --- that the linear inequality system given by Behrend and Knight~\cite{behrend2007higher} and by Striker~\cite{striker2007alternating, striker2009alternating} for describing the polytope of alternating sign matrices is totally dual integral (TDI), confirming a recent conjecture of Edmonds~\cite{Edmonds2024_talk2, Edmonds2024_talk1}.
  By relying on the polymatroidal approach, we derive a characterization for the existence of prefix-bounded matrices meeting lower and upper bounds on their entries.

  Furthermore, we point out that the constraint matrix of the linear system describing the convex hull of PBMs, in particular ASMs, is a network matrix.
  This implies that (a) standard network-flow techniques can be used to manage algorithmically optimization and structural results on PBMs obtained via g-polymatroids, (b) the linear system is actually box-TDI, and (c) the convex hull of PBMs admits a sharpened form of the integer Carath\'eodory property, in particular, the integer decomposition property.
  This latter feature makes it possible to confirm an extended form of an elegant conjecture of Brualdi and Dahl~\cite{brualdi2023multi} on the decomposability of a so-called $k$-regular alternating sign matrix as the sum of $k$ pattern-disjoint ASMs.\\



\end{abstract}

\section{Introduction}

A $(0,\pm 1)$-valued square matrix is called an {\bf alternating sign matrix} (ASM) if (a) the total sum of the entries in each row and in each column is $+1$, and (b) the $+1$ and $-1$ entries in each row and in each column alternately follow each other (implying that both the first and the last non-zero entries in each row and column are $+1$).
An immediate observation shows that a $(0,\pm 1)$-valued square matrix is an ASM if and only if (a) holds, and (b') the sum of the entries in every starting segment (which is called a prefix) of each row and each column is at least $0$ and at most $+1$.

The first appearance of this concept seems to be the paper of Mills et al.~\cite{mills1982proof} from 1982, where the authors already posed the ASM Conjecture concerning the number of ASMs of given size --- a conjecture that has since been confirmed~\cite{zeilberger1996}.
The study of ASMs is motivated by their connection with models in statistical mechanics, most notably the six-vertex model, also known as the square ice model, see the book of Bressoud~\cite{bressoud1999proofs}.

Since the first appearance of ASMs, a huge literature of the topic has been developed.
At the beginning, researchers concentrated mainly on enumerative problems concerning ASMs.
In 1999, Bressoud and Propp published a survey~\cite{bressoud1999alternating} detailing early results in this direction.
They write:

\medskip
{\it ``The ASM Conjecture has served to cross-fertilize the various modern offspring of classical invariant theory, drawing attention to connections no one had recognized.
  The study of alternating sign matrices should continue to bear fruit for many years to come and to tantalize us with fruit that is just beyond our reach.''}
\medskip


In 2007, Behrend and Knight~\cite{behrend2007higher} and Striker~\cite{striker2007alternating, striker2009alternating} opened a new door to investigations of ASMs by providing a polyhedral description of the polytope (convex hull) of ASMs.
Since then, various extensions of the notion of ASMs have been introduced and investigated, concerning, among others, polyhedral descriptions and existence (or feasibility) problems in which the goal is to explore situations where there exist (extended) ASMs with specified restrictions and prescriptions on their entries.
In this direction, we should emphasize the rich and fertilizing contribution of R.\ Brualdi and his co-authors.

For example, in a recent paper, Brualdi and Dahl~\cite{brualdi2024frobenius} solved a special case of the problem of characterizing $(0,\pm 1)$-valued square matrices which can be made an ASM by replacing some of their $\pm 1$ entries with $0$.
They pose the general problem of characterizing such $(0,\pm 1)$-valued square matrices, and write:

\medskip
{\it ``Basic Question:  Given an $n\times n$ $(0,\pm 1)$-matrix, when does it contain an $n\times n$ ASM?''}
\medskip

\noindent Here a matrix $X$ ``contains'' $A$, or in their more formal terms, $A$ is a subordinate of $X$ if $A$ can be obtained from $X$ by replacing some non-zero entries of $X$ with $0$.
Brualdi and Dahl solved the problem when $X$ belongs to a special class of matrices, see Theorem~4.1 in their paper, and remarked for the general case --- when $X$ is an arbitrary $(0,\pm 1)$-matrix --- that:

\medskip
{\it ``In this generality, there is probably no simple answer\dots''}
\medskip

In the present paper, we do provide such an answer in Corollary~\ref{cor:alarendelt}, which states, roughly, that the matrix $X$ has an alternating sign matrix subordinate if and only if it is not possible to cover the $+1$ entries of $X$ using fewer than $n$ separated vertical and horizontal segments (i.e., a subset of subsequent entries) in such a way that each $-1$ entry is covered by at most one segment.
Here we call a family of segments separated if any two segments of the same row or column are disjoint and have at least one entry between them (for the precise definition, see Section~\ref{sec:pb-g-poli}).
This result follows from the characterization in Theorem~\ref{thm:AS-feas2}, which establishes conditions for the existence of alternating sign matrices with certain entries fixed to specified values and others constrained by prescribed signs.
In fact, Theorem~\ref{thm:AS-feas2} will be formulated and proved in a broader context that extends beyond alternating sign matrices, see Theorem~\ref{thm:general1}.

What is the novelty that gave rise to this new result?
In 2024 and also in 2025, Jack Edmonds gave two conference talks entitled ``Ray and me''~\cite{Edmonds2024_talk2, Edmonds2024_talk1} on the occasion of the 100th birthday of D.R.\ Fulkerson.
In his slides, as well as in a circular email sent to several experts, he drew the attention of operations researchers by indicating that classical concepts and results of combinatorial optimization could be applicable to the study of alternating sign matrices.
For example, he formulated the following:

\medskip
{\it ``Ray Fulkerson's Birthday Conjecture (RFB): The Behrend-Knight-Striker linear inequality system is TDI.''  }
\medskip

\noindent Here TDI is the standard abbreviation of total dual integrality, which was introduced by Edmonds and Giles~\cite{edmonds1977min} and became a fundamental concept in combinatorial optimization, see also~\cite{schrijver2003combinatorial}.
Edmonds also remarked that the problem of finding a minimum cost ASM had not previously been investigated, and indicated that he found such an algorithm.
With this conjecture and remark, Edmonds opened a new door from combinatorial optimization to the hall of ASMs.
It was his initiative that triggered our present work, which is an attempt to enter this door.

The first step on this road is the natural observation that the polyhedral description of ASMs, given by Behrend, Knight, and Striker, can be viewed as a linear inequality system defined by the union of two special laminar families of sets --- for definitions, see the Preliminaries below.
We consider two different approaches to handling optimization and feasibility problems constrained by the union of two such laminar families.
These will be referred to as the polymatroidal view and the network-flow view.

\medskip

The fundamental concept of polymatroids was introduced and studied by Edmonds~\cite{edmonds1970submodular}, and later extended to generalized polymatroids (g-polymatroids, for short) in~\cite{frank1984generalized}.
The polymatroidal view of ASMs stems from the known fact that a polyhedron defined by arbitrary integer-valued lower and upper bounds on the members of a laminar family is an integral g-polymatroid.
This implies that the polytope of ASMs is the intersection of two g-polymatroids.

The realization that arbitrary lower and upper bounds on prefixes can be accommodated naturally led us to introduce the notion of prefix-bounded matrices (PBMs), a generalization of several previously studied extensions of ASMs.
Roughly speaking, a PBM is an integer-valued $m \times n$ matrix in which the sum of entries in each horizontal and vertical prefix lies between specified lower and upper bounds, see Section~\ref{sec:PBM}.
The properties (a) and (b') mentioned earlier in the definition of ASMs are special cases of such prefix bounds.

The convex hull of PBMs, although not necessarily bounded like the polytope of ASMs, is also the intersection of two (integral) g-polymatroids.
This allows us to apply notions and results from the theory of g-polymatroids, see Section~\ref{sec:generalized-polymatroids} for an overview.
For example, we obtain in this way a natural necessary and sufficient condition for the existence of a PBM satisfying lower and upper bounds on its entries, see Theorem~\ref{thm:general1}.
The ``Basic Question'' of Brualdi and Dahl cited earlier turns out to be a very special case of this theorem.

We emphasize that the polymatroidal view played a central role in formulating the correct condition in Theorem~\ref{thm:general1}.
Although our first proof uses the g-polymatroid intersection theorem, the statement of the theorem itself does not refer to g-polymatroids.
Since there are purely combinatorial polynomial-time algorithms for deciding whether the intersection of two g-polymatroids is empty, the feasibility problems of PBMs can be managed algorithmically.
Nevertheless, it is natural to seek algorithms that avoid polymatroidal concepts and rely instead on simpler tools such as standard network-flow algorithms.

\medskip

The network-flow view of ASMs stems from the following earlier results.
Edmonds and Giles~\cite{edmonds1977min} proved that a laminar family can be represented by an arborescence, see Theorem~1.4.1 in~\cite{frank2011connections} for a reformulation.
This, in turn, implies Theorem~4.2.8 in~\cite{frank2011connections}, which asserts that the transpose of the incidence matrix of the union of two laminar families is a network matrix.
An earlier result of Tutte~\cite{tutte1965lectures} had already shown that every network matrix is TU, implying that the system $Ax\leq b$ is TDI for any integer-valued bounding vector $b$ --- which also follows from Theorem (39) in~\cite{edmonds1970submodular} due to Edmonds.
This immediately settles Edmonds' RFB conjecture not only for ASMs but also for PBMs, see Corollary~\ref{cor:box-TDI}.
Combining the TU property of $A$ with a simple new observation stated in Proposition~\ref{prop:strong-Cara}, we also prove a recent conjecture of Brualdi and Dahl~\cite{brualdi2023multi}, see Corollary~\ref{cor:B+D-conj}.

An important consequence of these facts is that various optimization and feasibility problems concerning PBMs can be solved using standard network-flow techniques.
Namely, with the help of a suitable construction, Theorem~\ref{thm:general1} --- first proved via g-polymatroids --- will also be derived from Hoffman's classical theorem on the existence of feasible circulations without relying on g-polymatroids.

\medskip

We also introduce a further extension of ASMs that lies outside the PBM framework.
Within each row and column, we impose two types of constraints: lower and upper bounds on the total sum of the entries in that row or column, and uniform lower and upper bounds (depending on the row or column) on the sum of entries in every segment.
In addition, lower and upper bounds are prescribed for the entries.
As noted in Remark~\ref{rem:uniform_segment}, the convex hull of such integer matrices is again the intersection of two g-polymatroids.
ASMs arise as the special case where the uniform upper bound for segments is $+1$, the lower bound for full rows and columns is $+1$, and the lower bound for other segments is $-1$.
This model behaves quite differently from PBMs: whereas the prefix-bounded setting admits a network-flow formulation, allowing feasibility and optimization problems to be solved by standard algorithms, the segment-bounded framework appears to lack such a reduction, and, to the best of our knowledge, the only available method relies on g-polymatroids.

\medskip

The strength of our approach is highlighted by the fact that feasibility problems become NP-complete for certain natural extensions of PBMs.
In particular, NP-completeness arises when bounds are imposed not only on prefixes of rows and columns but also on suffixes --- that is, on the sums of entries in terminal segments of rows or columns --- along with entry bounds.
Furthermore, the segment-bounded problem under entry bounds becomes NP-complete when, within each row and column, uniform lower and upper bounds (depending on the row or column) are imposed on the sum of the entries in the segments of a single length specified in the input --- rather than on all segments as above.

\medskip

The rest of the paper is organized as follows.
Section~\ref{sec:prelim} introduces the notation and preliminary concepts used throughout the paper.
Section~\ref{sec:generalized-polymatroids} reviews the necessary background on generalized polymatroids, which underpins our structural results.
In Section~\ref{sec:special-g-polymatroids}, we introduce the notions of laminar and prefix-bounded g-polymatroids, which serve as the main tools in our study of PBMs.
Section~\ref{sec:intersection} establishes intersection theorems for two g-polymatroids under additional box and plank constraints.
Section~\ref{sec:PBM} formally introduces prefix-bounded matrices (PBMs) and derives their fundamental structural properties.
In Section~\ref{sec:special-PBM}, we examine notable subclasses of PBMs --- including ASMs --- and explore their connections to the six-vertex model.
Finally, Section~\ref{sec:circulation} relates PBMs to a feasible circulation problem and provides another proof of Theorem~\ref{thm:general1}, which is self-contained, algorithmic, and does not rely on g-polymatroids.

\subsection{Preliminaries}\label{sec:prelim}

For an integer $j\geq 1$, $[j]$ denotes the set $\{1,\dots ,j\}$.
For integers $i\leq j$, $[i,j] \coloneqq \{i,\dots ,j\}$.
For integers $m\geq 1, n\geq 1$, let $S_{m,n} \coloneqq [m] \times [n]$.
$S_{m,n}$ may also be viewed as a two-dimensional array.
A matrix $A$ of size $m \times n$ may be viewed as a two-variable function on the set $S_{m,n}$.
Therefore, the entry of $A$ in row $i$ and column $j$ is denoted by $A(i,j)$.

For a function $f$ on ground set $S$ and for $X\subseteq S$, we use the notation
\[
  \widetilde f(X) \coloneqq \sum\ [f(s) : s\in X],
\]
where $\widetilde f(\emptyset) \coloneqq 0$.
In particular, for a set function $p$ and a family of sets $\cal I$, we have $\widetilde p({\mathcal I}) \coloneqq \sum [p(X) : X\in {\cal I}]$.
For set functions, we assume throughout that they are $0$ on the empty set.

The union $X\cup \{v\}$ of a set $X$ and a singleton $\{v\}$ is denoted by $X+v$.
The difference $\{v :  v\in X, v \notin Y\}$ of sets $X$ and $Y$ is denoted by $X-Y$.
The symmetric difference of sets $X$ and $Y$ is denoted by $X \ominus Y \coloneqq (X-Y)\cup (Y-X)$.
The complement $S-Z$ of a subset $Z\subseteq S$ is denoted by $\ol Z$.

Two subsets $X$ and $Y$ of a ground set $S$ are called {\bf intersecting} if $X\cap Y \neq \emptyset$ and {\bf properly intersecting} if $X-Y$, $Y-X$, and $X\cap Y$ are non-empty.
If, in addition, $S-(X\cup Y)$ is non-empty, then $X$ and $Y$ are {\bf crossing}.
A family of subsets is {\bf laminar} (resp., {\bf cross-free}) if it contains no two properly intersecting (crossing) members.
A {\bf sub-partition} $\{X_1,\dots ,X_k\}$ of ground set $S$ is a family of disjoint subsets of $S$.

Although the concept of a network matrix is well-known, see the books~\cite{frank2011connections, schrijver2003combinatorial}, we include its definition here for completeness. 

Let $T$ be the arc set of a (directed) spanning tree of a (weakly) connected digraph $D = (V,A)$, and let $N \coloneqq A-T$ be the set of non-tree arcs.
A {\bf network matrix} $M_D$ belonging to $T$ is a $(0,\pm 1)$-valued $\vert T\vert \times \vert N\vert$ matrix in which the rows (columns) correspond to the elements of $T$ \ ($N$).
The entry $M_D(e,f)$ belonging to row $e\in T$ and column $f\in N$ is non-zero if and only if $e$ lies in the unique circuit $C$ of $T + f$ (the fundamental circuit of $f$); in that case, $M_D(e,f) = +1$ if $e$ and $f$ are oriented in the same direction in $C$, and $-1$ otherwise.
In particular, if every fundamental circuit happens to be a one-way circuit, then the network matrix is $(0,+1)$-valued.

It is known that if a linear cost function $cx$ on a polyhedron $Q$ (which is the intersection of a finite number of closed half-spaces, or in other words, the solution set of a system of linear inequalities) is bounded from below, then $\min \{cx :  x\in Q\}$ exists.
When $cx$ is not bounded from below over $Q$, we write $\min \{cx :  x\in Q\} = -\infty$, and similarly if $cx$ is not bounded from above, then $\max \{cx :  x\in Q\} = +\infty$.

A {\bf box} $T(f,g) \coloneqq \{x\in \R^S : f\leq x\leq g\}$ is a polyhedron defined by functions $f\leq g$ where $f$ may have $-\infty$ values, $g$ may have $+\infty$ values.
A {\bf plank} $K(\alpha,\beta) \coloneqq \{x\in \R^S : \alpha \leq \widetilde x(S) \leq \beta \}$ is a polyhedron defined by numbers $-\infty \leq \alpha \leq \beta \leq +\infty$.

For a polyhedron $Q \coloneqq \{x : Ax\leq b\}$ and a positive integer $k$, the polyhedron $\{x : Ax\leq kb\}$ is called the {\bf $\boldsymbol{k}$-elongation} (or just an elongation) of $Q$.
An integral polyhedron $Q$ is said to admit the {\bf integer decomposition property} if, for every positive integer $k$, every integral vector of the $k$-elongation of $Q$ can be expressed as the sum of $k$ integral vectors in $Q$.
Baum and Trotter~\cite{baum1978integer} proved that a polyhedron described by a TU matrix and an integral bounding vector has the integer decomposition property.

As a significant sharpening of this result, Gijswijt and Regts~\cite{gijswijt2012polyhedra} proved that the projection (along coordinate axes) of a polyhedron $Q \subseteq \R^S$ defined by a TU matrix and an integral bounding vector admits the {\bf integer Carath\'eodory property}, see Theorem~6 in their paper.
That is, for every positive integer $k$, any integral vector $z^*$ in the $k$-elongation of $Q$ can be expressed as a positive integer linear combination of affinely independent integral vectors in $Q$.
More precisely, there exist affinely independent integral vectors $z_1, \dots, z_t \in Q$ and positive integers $\mu_1, \dots, \mu_t$ such that $\mu_1 + \dots + \mu_t = k$ and $z^* = \mu_1 z_1 + \dots + \mu_t z_t$.
In particular, every integral vector $z^*$ in the $k$-elongation of $Q$ can be expressed as a positive integer linear combination of at most $|S|+1$ integral members of $Q$.

We say that in the expression $z^* = \mu_1z_1 + \cdots + \mu_tz_t$, the vectors $z_1, \dots ,z_t$ are {\bf sign-consistent with} $z^*$ if, for each $s \in S$, $z^*(s) \geq 0$ implies $z_i(s) \geq 0$ for each $i \in [t]$, and $z^*(s) \leq 0$ implies $z_i(s) \leq 0$ for each $i \in [t]$.
We close this section by showing that the theorem of Gijswijt and Regts has a self-refining nature, as captured in the following proposition.
\begin{proposition}\label{prop:strong-Cara}
  Let $A$ be a TU matrix and $b$ an integer vector such that the polyhedron $Q \coloneqq \{x: Ax\leq b\}\subseteq \R^S$ is non-empty.
  For a positive integer $k$, let $Q_k$ denote the $k$-elongation of $Q$, that is, $Q_k \coloneqq \{x:  Ax\leq kb \}\subseteq \R^S$, and let $z^*$ be an integral vector in $Q_k$.
  Then there exist affinely independent integral vectors $z_1,\dots ,z_t \in Q$ that are sign-consistent with $z^*$, and there are positive integer coefficients $\mu _1, \dots, \mu_t$ with $\mu_1 + \dots + \mu_t = k$ such that $z^* = \mu_1z_1 + \cdots + \mu_tz_t$.
\end{proposition}
\begin{proof}
  We show that the proposition follows directly from the theorem of Gijswijt and Regts.
  To see this, extend the original linear system $Ax\leq b$ by the following inequalities for all $s\in S$:
  \[
    \hbox{ $x(s)\geq 0$ \ if \ $z^*(s)\geq 0$, \quad \quad \quad $x(s)\leq 0$ \ if \ $z^*(s)\leq 0$.}
  \]

  The resulting system still has a totally unimodular constraint matrix, and it defines a new polyhedron $Q' \subseteq Q$.
  Let $Q'_k$ denote the $k$-elongation of $Q'$.
  It immediately follows that $z^*\in Q'_k$, thus, by applying the original theorem of Gijswijt and Regts, we obtain affinely independent integral vectors $z_1, \dots ,z_t \in Q'$, and positive integer coefficients $\mu _1, \dots, \mu_t$ with $\mu_1 + \dots + \mu_t = k$ such that $z^* = \mu _1z_1 + \cdots + \mu _tz_t$.
  The definition of $Q'$ ensures that the vectors $z_1, \dots, z_t$ are sign-consistent with $z^*$, as claimed.
\end{proof}

Proposition~\ref{prop:strong-Cara} may be considered as a sign-consistent sharpening of the theorem of Gijswijt and Regts.
The proposition will be used to derive Corollaries~\ref{cor:integer-Cara}~and~\ref{cor:int-decomp-spec} which result in the proof of a conjecture of Brualdi and Dahl~\cite{brualdi2023multi}.

\section{Generalized polymatroids}\label{sec:generalized-polymatroids}

As already mentioned in the Introduction, the main driving force behind the present work is the simple yet fundamental observation that the convex hull of alternating sign matrices is the intersection of two (rather special) generalized polymatroids.
This is why we outline some important properties of g-polymatroids in the following two sections.
As a natural extension of polymatroids, g-polymatroids were introduced and studied by Frank~\cite{frank1984generalized}.
The books of Frank~\cite{frank2011connections} and Schrijver~\cite{schrijver2003combinatorial} are rich sources for further information.
This polymatroidal background helps formulating the right theorems, which is useful even in cases when it is possible to work out proofs not relying on g-polymatroids.
We also hope that the polymatroidal view shall open new doors toward further investigations.

In what follows, we are given two integer-valued set functions $p$ and $b$, defined on a ground set $S$, as well as two integer-valued functions $f$ and $g$ on $S$.
Here, $f$ and $p$ serve as lower bounds, and their values are permitted to be $-\infty$.
Analogously, the upper bounds $g$ and $b$ are allowed to take $+\infty$ values.

A set function $b$ on ground set $S$ is {\bf fully submodular} or just {\bf submodular} if the submodular inequality
\[
  b(X) + b(Y) \geq b(X\cap Y) + b(X\cup Y)
\]
holds for every $X,Y \subseteq S$.
When the submodular inequality is required only for intersecting (respectively, crossing) pairs $X,Y$, we speak of an {\bf intersecting (crossing) submodular} function.
A set function $p$ is {\bf fully (intersecting, crossing) supermodular} if $-p$ is fully (intersecting, crossing) submodular.

A pair $(p, b)$ of set functions on a ground set $S$ is called {\bf strongly paramodular} (or, simply, a {\bf strong pair}) if $p$ is fully supermodular, $b$ is fully submodular, and the {\bf cross-inequality}
\[
  b(X) - p(Y) \geq b(X - Y) - p(Y - X)
\]
holds for all $X, Y \subseteq S$.

For a strong pair $(p,b)$, the polyhedron $Q(p,b) \coloneqq \{x : p(Z) \leq \widetilde x(Z) \leq b(Z)$ for every $Z \subseteq S \}$ (or, more concisely, $\{x : p\leq \widetilde x\leq b\}$) is called a {\bf g-polymatroid}, and we say that $Q(p,b)$ is {\bf bordered} by $(p,b)$.
The empty set is also considered to be a g-polymatroid.
It is known that both a face and a projection of a g-polymatroid are also g-polymatroids.
A deeper theorem asserts that the (Minkowski-) sum of g-polymatroids $Q(p_i,b_i)$ \ ($i = 1,\dots ,k)$ bordered by integer-valued strong pairs is also a g-polymatroid whose unique bordering strong pair is $(\sum p_i, \sum b_i)$, see Theorem~14.2.15~in~\cite{frank2011connections}.

Let $b$ be a fully submodular function for which $b(S)$ is finite.
The polyhedron $B(b) \coloneqq \{x : \widetilde x(Z)\leq b(Z)$ for each $Z\subset S$ and $\widetilde x(S) = b(S)\}$ is called a {\bf base-polyhedron}.
When $b(S) = 0$, we speak of a {\bf 0-base-polyhedron}.
The empty set is also considered to be a base-polyhedron.
When $p$ is a fully supermodular function for which $p(S)$ is finite, the polyhedron $B'(p) \coloneqq \{x : \widetilde x(Z)\geq p(Z)$ for every $Z\subset S$ and $\widetilde x(S) = p(S)\}$ is a base-polyhedron, namely, $B'(p) = B(b)$ where $b(X) \coloneqq p(S)-p(S-X)$ for every $X\subseteq S$.
\begin{claim}\label{cl:pSbS}
  Let $(p^*,b^*)$ be a strong pair for which $p^*(S) = b^*(S)$.
  Then $p^*(X) = b^*(S)-b^*(S-X)$ for every subset $X\subseteq S$, and the g-polymatroid $Q(p^*,b^*)$ is the base-polyhedron $B(b^*) = B'(p^*)$.
\end{claim}
\begin{proof}
  By applying the cross-inequality first to sets $S-X, S$ and then to sets $S,X$, we get that $p^*(S) - b^*(S-X) \leq p^*(X)$ and $b^*(S) - p^*(X) \geq b^*(S-X)$, and hence $p^*(S) - b^*(S-X) \leq p^*(X) \leq b^*(S) - b^*(S-X)$.
  But here we must have equality because $p^*(S) = b^*(S)$.
  Therefore, $p^*(X) = b^*(S) - b^*(S-X)$, and we get similarly that $b^*(X) = p^*(S) - p^*(S-X)$, from which $Q(p^*,b^*) = B(b^*) = B'(p^*)$ follows.
\end{proof}

The following result is also known from the literature, see Proposition~II.2.2 in~\cite{frank1988generalized}.
\begin{proposition}
  Let $(p^*, b^*)$ be an integer-valued strong pair.
  Then the g-polymatroid $Q(p^*,b^*)$ is an integral polyhedron which is never empty.
  \FBOX
\end{proposition}

A pair $(p,b)$ of set functions is called {\bf weakly paramodular} (or, simply, a {\bf weak pair}) if $p$ is intersecting supermodular, $b$ is intersecting submodular, and the cross-inequality holds for properly intersecting subsets $X,Y \subseteq S$.
By Propositions~II.2.5 and~II.2.6 in~\cite{frank1988generalized}, we have
\begin{proposition}\label{prop:j17-2.6}
  {\bf (A)} \ The polyhedron $Q(p,b)$ bordered by a weak pair $(p,b)$ is a (possibly empty) g-polymatroid, which is a base-polyhedron when $p(S) = b(S)$.

  \medskip

  \noindent {\bf (B)} \ The polyhedron $Q(p,b)$ is non-empty if and only if $p(\cup Z_i) \leq \sum b(Z_i)$ and $b(\cup Z_i) \geq \sum p(Z_i)$ hold for every sub-partition $\{Z_1,\dots ,Z_q\}$ of the ground set $S$.
  \FBOX
\end{proposition}

Part (B) of Proposition~\ref{prop:j17-2.6} is, in fact, equivalent to a theorem of Fujishige~\cite{fujishige1984structures}, which characterizes the non-emptiness of a base-polyhedron bordered by a crossing submodular function.

The following four propositions on g-polymatroids are useful results from the literature, see the paper~\cite{frank1988generalized} or the book~\cite{frank2011connections}.

\begin{proposition}\label{prop:hatarpar}
  Let $Q$ be a non-empty g-polymatroid.
  There exists a unique strong pair $(p^*,b^*)$ for which $Q = Q(p^*,b^*)$, namely,
  \[
    p^*(Z) = \min \{ \widetilde x(Z) : x\in Q\} \ \ \ \hbox{and}\ \ \ b^*(Z) = \max \{ \widetilde x(Z) : x\in Q\}.
  \]
  \FBOX
\end{proposition}

\begin{proposition}
  Suppose that a g-polymatroid $Q \coloneqq Q(p,b)$ bordered by a weak pair $(p,b)$ is non-empty and that $p(S) = b(S)$.
  Let $(p^*,b^*)$ denote the unique strong pair bordering $Q$.
  Then $Q = B(b^*)$.
  \FBOX
\end{proposition}

For a weak pair $(p,b)$, it is easily seen that if we increase the $p$-values on some singletons and/or decrease the $b$-values on some singletons, then the modified pair $(p^+,b^-)$ is also weakly paramodular.
It is also true that we get a weak pair when $p(S)$ is increased and/or $b(S)$ is decreased.
These imply the following.

\begin{proposition}\label{prop:box_plank_intersection}
  The intersection of a g-polymatroid $Q$ with a box $T(f,g) = \{x : f\leq x\leq g\}$ is a g-polymatroid, and so is the intersection of $Q$ with a plank $K(\alpha,\beta ) = \{x : \alpha \leq \widetilde x(S) \leq \beta \}$.
  \FBOX
\end{proposition}

Proposition~II.2.11 in~\cite{frank1988generalized} is as follows.
\begin{proposition}\label{prop:j17-2.11}
  Suppose that the intersection of a g-polymatroid $Q \coloneqq Q(p^*,b^*)$ bordered by a strong pair $(p^*,b^*)$ with a box $T \coloneqq T(f,g) = \{x : f\leq x\leq g\}$ is non-empty.

  \medskip

  \noindent {\bf (A)} \ The unique strong pair $(p^T,b^T)$ bordering the g-polymatroid $Q\cap T$ is
  \[
    \begin{cases}
      & p^T (Z) = \max \{p^*(X) + \widetilde f(Z-X) - \widetilde g(X-Z) : X\subseteq S \},\\
      & b^T (Z) = \min \{b^*(X) + \widetilde g(Z-X) - \widetilde f(X-Z) : X\subseteq S \}.
    \end{cases}
  \]

  \noindent {\bf (B)} \ In the special case when $g\equiv +\infty$, the formula in (A) reduces to
  \[
    \begin{cases}
      & p^T (Z) = \max \{p^*(X) + \widetilde f(Z-X) : X\subseteq Z \},\\
      & b^T (Z) = \min \{b^*(X) - \widetilde f(X-Z) : X \supseteq Z \}.
    \end{cases}
  \]

  \noindent {\bf (C)} \ In the special case when $f\equiv -\infty$, the formula in (A) reduces to
  \[
    \begin{cases}
      & p^T (Z) = \max \{p^*(X) - \widetilde g(X-Z) : X \supseteq Z \},\\
      & b^T (Z) = \min \{b^*(X) + \widetilde g(Z-X) : X \subseteq Z \}.
    \end{cases}
  \]
  \FBOX
\end{proposition}





\section{Special g-polymatroids}\label{sec:special-g-polymatroids}

\subsection{Laminar g-polymatroids }\label{sec:lam.g-poli}

A special case of g-polymatroids bounded by a weak pair is the following.
Let $\cal L$ be a laminar family of non-empty subsets of ground set $S$, and suppose that the bordering set functions $p$ and $b$ may be finite only on the members of $\cal L$.
In this case, the intersecting sub- and supermodularity as well as the cross-inequality requires nothing since a laminar family has no two properly intersecting members.
We refer to such a g-polymatroid $Q(p,b)$ as a {\bf laminar g-polymatroid}.
Since adding singletons and/or the ground set $S$ to $\cal L$ preserves laminarity, the intersection of a laminar g-polymatroid with a box or a plank is again a laminar g-polymatroid.
Similarly, both faces and projections of a laminar g-polymatroid remain laminar g-polymatroids.

The following basic observation was formulated by Edmonds~\cite{edmonds1970submodular} as Theorem~(39) in his paper, see also Theorem~41.11 in the book of Schrijver~\cite{schrijver2003combinatorial}. 
\begin{lemma}[Edmonds]\label{lem:edmonds-lemma}
  Let ${\cal L}'$ be the union of two laminar families on a ground set $S$.
  Then the transpose of the incidence matrix of ${\cal L}'$ is totally unimodular, in which the rows correspond to the members of the family and the columns correspond to the ground set.
  \FBOX
\end{lemma}

\begin{remark}\label{rem:2-lamin}
  It is easy to show that the transpose of the incidence matrix of ${\cal L}'$ is a network matrix, see Theorem~4.2.8 in~\cite{frank2011connections}.
  A classical theorem of Tutte~\cite{tutte1965lectures} states that every network matrix is TU, implying Edmonds' lemma, see also Theorem~4.2.5 in~\cite{frank2011connections}.
  It is also known that a polyhedron described by an arbitrary network matrix is the projection of a feasible circulation polyhedron, see Theorem~4.3.2 in~\cite{frank2011connections}.
  This implies that a linear program described by a network matrix (or its transpose) can be formulated as a cheapest circulation problem (or its dual), and hence it can be solved via standard network-flow algorithms in strongly polynomial time.
  $\bullet$
\end{remark}

It is known that a linear inequality system described by a TU matrix and an integral bounding vector is box-TDI.
Since the intersection of a laminar g-polymatroid with a box and a plank is also a laminar g-polymatroid, the following theorem follows from Lemma~\ref{lem:edmonds-lemma} and Remark~\ref{rem:2-lamin}.
\begin{theorem}\label{thm:box-tdi1}
  Let ${\cal L}_1$ and ${\cal L}_2$ be two laminar families on ground set $S$.
  For $i = 1,2$, let $p_i \leq b_i$ be integer-valued functions on ${\cal L}_i$.
  Moreover, let $f\leq g$ be lower and upper bound functions on $S$, and let $\alpha \leq \beta$ be two integers.

  \medskip

  \noindent {\bf (A)} \ The linear inequality system
  \begin{equation}\label{eq:2-lamin-leir}
    \begin{cases}
      & p_1(Z) \leq \widetilde x(Z) \leq b_1(Z) \hbox{ for every } Z \in {\mathcal{L}}_1\\
      & p_2(Z) \leq \widetilde x(Z) \leq b_2(Z) \hbox{ for every } Z \in {\mathcal{L}}_2\\
      & f\leq x\leq g\\
      & \alpha \leq \widetilde x(S) \leq \beta
    \end{cases}
  \end{equation}
  is box-TDI.
  \medskip

  \noindent {\bf (B)} \ For a cost function $c$, a cheapest (integral) element of the polyhedron described by~\eqref{eq:2-lamin-leir} can be computed with the help of an algorithm for computing a cheapest feasible circulation.
  \FBOX
\end{theorem}

\subsection{Prefix-bounded g-polymatroids }\label{sec:pb-g-poli}

Let $U \coloneqq \{u_1,\dots ,u_t\}$ be a totally ordered ground set, and let $u_0$ be an additional element not contained by $U$.
A subset $U_{[h,k]} \coloneqq \{u_h,\dots ,u_k\}$ \ $(1\leq h\leq k\leq t)$ of $U$ consisting of subsequent elements of $U$ is called a {\bf segment}, while the special segment $U_{[k]} \coloneqq U_{[1,k]}$ is a {\bf prefix} of $U$.
An arbitrary subset $Z\subseteq U$ can uniquely be decomposed into disjoint maximal segments $Z_1,\dots ,Z_q$, which are {\bf separated} in the sense that they do not ``touch'' each other, that is, they are separated by some elements of $U-Z$.
We refer to these subsets $Z_i$ as the {\bf segments in $\boldsymbol Z$}.
Figure~\ref{fig:prefix-segments} illustrates these definitions.

\medskip
\begin{figure}[h!]
  \centering
  \begin{tikzpicture}[
    node/.style={font=\small},
    every path/.style={thick},
    brace/.style={decorate, decoration={brace, amplitude=8pt}},
    brace_m/.style={decorate, decoration={brace, mirror, amplitude=8pt}},
    arrow/.style={-{Stealth[length=5pt]}, thick}
    ]

    \node[node] (u1) at (0,0) {$u_1$};
    \node[node] (u2) at (1.5,0) {$u_2$};
    \node[node] (u3) at (3.0,0) {$u_3$};
    \node[node] (u4) at (4.5,0) {$u_4$};
    \node[node] (u5) at (6.0,0) {$u_5$};
    \node[node] (u6) at (7.5,0) {$u_6$};
    \node[node] (u7) at (9,0) {$u_7$};

    \draw[blue!90,opacity=0.6,line width=18pt,line cap=round] (u1.west) -- (u3.east) node[text=black, opacity=1, midway, yshift=-20pt] {Prefix $U_{[3]}$};;
    \draw[blue!90,opacity=0.6,line width=18pt,line cap=round] (u5.west) -- (u6.east) node[text=black, opacity=1, midway, yshift=-20pt] {Segment $U_{[5,6]}$};
  \end{tikzpicture}
  \caption{Illustration of the prefix $U_{[3]} = \{u_1, u_2, u_3\}$ and the segment $U_{[5,6]} = \{u_5, u_6\}$ in a totally ordered ground set $U = \{u_1, \dots, u_7\}$.
    Together they form the unique decomposition of the set $\{u_1,u_2,u_3,u_5,u_6\}$ into separated maximal segments.}\label{fig:prefix-segments}
\end{figure}

Let $\varphi \leq \gamma$ be two integer-valued functions on $U+u_0$, for which $\varphi (u_0) = \gamma (u_0) = 0$.
Here $\varphi$ may take $-\infty$ values and $\gamma$ may take $+\infty$ values.
Define two set functions $p_\varphi$ and $b_\gamma$ on $U$ as follows.
For a prefix $Z \coloneqq U_{[k]}$, let $p_\varphi (Z) \coloneqq \varphi (u_k)$ and $b_\gamma (Z) \coloneqq \gamma (u_k)$ \ ($k = 1,\dots ,t$), and, for all other non-empty subsets $Z$ of $U$, let $p_\varphi (Z) \coloneqq -\infty$ and $b_\gamma (Z) \coloneqq +\infty$.
Then $(p_\varphi,b_\gamma)$ is a weak pair, where $p_\varphi$ is fully supermodular, $b_\gamma$ is fully submodular, and $p_\varphi \leq b_\gamma$.
By the discrete separation theorem~\cite{frank1982algorithm}, there always exists a modular function between $p_\varphi$ and $b_\gamma$, thus the g-polymatroid $Q(p_\varphi,b_\gamma)$ is non-empty.
We call such a g-polymatroid an {\bf elementary prefix-bounded g-polymatroid}, while the direct sum of such g-polymatroids is a {\bf prefix-bounded g-polymatroid}.
It follows from the definition that a prefix-bounded g-polymatroid is laminar; consequently, its intersection with a box and a plank is also laminar.
Our objective is to identify the unique strong pair that borders it.
To do so, it suffices to consider only elementary prefix-bounded g-polymatroids.
\medskip

Define two set functions $p^*$ and $b^*$ on $U$ as follows.
For a segment $Z \coloneqq U_{[h,k]}$ \ ($1\leq h\leq k\leq t$), let
\begin{equation}\label{eq:p'b'def}
  \begin{cases}
    & p^*(Z) \coloneqq \varphi (u_k) - \gamma (u_{h-1}),\\
    & b^*(Z) \coloneqq \gamma (u_k) - \varphi (u_{h-1}),
  \end{cases}
\end{equation}
and for a non-empty subset $Z\subseteq U$ with segments $Z_1,\dots ,Z_q$, let
\begin{equation}\label{eq:addit}
  \begin{cases}
    & p^*(Z) \coloneqq \sum\ [p^*(Z_r) : r = 1,\dots, q],\\
    & b^*(Z) \coloneqq \sum\ [b^*(Z_r) : r = 1,\dots, q].
  \end{cases}
\end{equation}

\begin{lemma}\label{lem:p*b*}
  The pair $(p^*,b^*)$ of set functions defined in~\eqref{eq:addit} is the unique strong pair bordering the elementary prefix-bounded g-polymatroid $Q(p_\varphi, b_\gamma) \subseteq \R^U$.
  If $\varphi (u_t) = \gamma (u_t)$, then $Q(p_\varphi, b_\gamma)$ is the base-polyhedron for which $p^*(U) = \varphi (u_t) = \gamma (u_t) = b^*(U)$ and $p^*(Z) = b^*(U) - b^*(U-Z)$.
\end{lemma}
\begin{proof}
  We present the details of the proof only for the special case when both $\varphi$ and $\gamma$ are finite-valued.
  The general case follows by similar arguments, though it involves additional technical considerations, see Remark~\ref{rem:infdet}.

  By Proposition~\ref{prop:hatarpar}, what we have to prove for every subset $Z\subseteq U$ is that
  \[
    \begin{cases}
      & p^*(Z) = \min \{\widetilde x(Z) : x\in Q(p_\varphi, b_\gamma)\},\\
      & b^*(Z) = \max \{\widetilde x(Z) : x\in Q(p_\varphi, b_\gamma)\}.
    \end{cases}
  \]

  We prove only the first equality, as the second one follows analogously.
  Let us prove first that $\min \{\widetilde x(Z) : x\in Q(p_\varphi, b_\gamma)\} \geq p^*(Z)$, that is, $\widetilde x(Z) \geq p^*(Z)$ holds for an arbitrary element $x\in Q(p_\varphi, b_\gamma)$.

  Consider first the case when $Z$ itself is a segment $U_{[h,k]}$ \ ($1\leq h\leq k \leq t$).
  If $h = 1$, then for the prefix $U_{[k]} = U_{[1,k]}$, we have
  \[
    \widetilde x(Z) =  \widetilde x(U_{[k]})\geq p_\varphi (U_{[k]}) = \varphi (u_k) = p^*(U_{[k]}) = p^*(Z).
  \]

  If $h\geq 2$, then for a non-prefix segment $U_{[h,k]}$, we have
  \begin{align}\label{eq:x(z)}
    \begin{split}
      \widetilde x(Z) & = \widetilde x(U_{[h,k]}) = \widetilde x(U_{[k]}) - \widetilde x(U_{[h-1]}) \geq  p_\varphi (U_{[k]}) - b_\gamma (U_{[h-1]})
      \\ & = \varphi (u_k) - \gamma (u_{h-1}) = p^*(U_{[h,k]}) = p^*(Z).
    \end{split}
  \end{align}

  In the general case, let $Z_1,\dots ,Z_q$ be the maximal segments in $Z$.
  By applying~\eqref{eq:x(z)} to $Z_r$ \ ($r = 1,\dots ,q$) in place of $Z$, we obtain the following from the definition of $p^*$ given in~\eqref{eq:addit}:
  \[
    \widetilde x(Z) = \sum\ [\widetilde x(Z_r) : r = 1,\dots ,q] \geq \sum\ [p^*(Z_r) : r = 1,\dots ,q] = p^*(Z).
  \]

  To see the reverse inequality $\min \{\widetilde x(Z) : x\in Q(p_\varphi, b_\gamma)\} \leq p^*(Z)$, we have to find an element $x_Z\in Q(p_\varphi, b_\gamma)$ for which $\widetilde x_Z(Z) = p^*(Z)$.
  To this end, define function $\mu _Z : U + u_0\rightarrow \Z$ as follows.
  \[
    \mu _Z(u) \coloneqq
    \begin{cases}
      \varphi (u) & \ \ \hbox{if}\ \ \ u\in Z,
      \\ \gamma (u) & \ \ \hbox{if}\ \ \ u\in U-Z,
      \\ 0 & \ \ \hbox{if}\ \ \ u = u_0.
    \end{cases}
  \]

  Obviously, $\varphi (u) \leq \mu _Z(u) \leq \gamma (u)$ holds for every element $u\in U$.
  Define $x_Z : U\rightarrow \Z$ as
  \[
    x_Z(u_j) \coloneqq \mu _Z(u_j)- \mu _Z(u_{j-1}) \ \ (j = 1,\dots ,t).
  \]
  For every segment $U_{[h,k]}$ \ ($1\leq h\leq k\leq t$), we have the following.
  \[
    \widetilde x_Z(U_{[h,k]}) = \mu _Z(u_k) - \mu _Z(u_{h-1}),
  \]
  and in particular
  \[
    \widetilde x_Z(U_{[k]}) = \mu _Z(u_k)
  \]
  holds for every prefix $U_{[k]} = U_{[1,k]}$.
  Therefore,
  \[
    p_\varphi (U_{[k]}) = \varphi (u_k) \leq \widetilde x_Z(U_{[k]}) \leq \gamma (u_k) = p_\gamma (U_{[k]}),
  \]
  that is, $x_Z\in Q(p_\varphi, b_\gamma)$.

  We claim that $\widetilde x_Z(Z) = p^*(Z)$.
  By~\eqref{eq:addit}, it suffices to see this only for the segments in $Z$.
  Let $U_{[h,k]}$ be a segment of $Z$.
  Then
  \[
    \widetilde x_Z(U_{[h,k]}) = \mu _Z(u_k) - \mu _Z(u_{h-1}) = \varphi (u_k) - \gamma (u_{h-1}) = p^*(U_{[h,k]}),
  \]
  and thus the proof of the first part of the lemma is complete.

  The second half of the lemma follows directly from Claim~\ref{cl:pSbS}.
\end{proof}

\begin{remark}\label{rem:infdet}
  The proof of Lemma~\ref{lem:p*b*} in the general case, when $\varphi$ may take $-\infty$ values and $\gamma$ may take $+\infty$ values, can roughly be done along the following line.
  Let $K$ be a sufficiently large positive integer.
  For every element $u\in U$ with $\varphi (u) = -\infty$, let $\varphi (u) \coloneqq -K$, and for every element $u$ with $\gamma (u) = +\infty$, let $\gamma (u) \coloneqq K$.
  The general case of the lemma follows by applying its finite-valued case to these modified (finite-valued) $\varphi$ and $\gamma$.
  $\bullet$
\end{remark}

\begin{remark}
  A special class of elementary prefix-bounded g-polymatroids is the following.
  Let $\gamma (u_i) \coloneqq +1$ for $i\in [t]$, and let $\varphi (u_i) \coloneqq 0$ for $i\in [t-1]$ and $\varphi (u_t) \coloneqq +1$.
  Then the prefix-bounded g-polymatroid defined by these functions is a base-polyhedron whose integral elements are exactly the $(0,\pm 1)$-valued, not-identically-zero vectors for which the first and last non-zero components are $+1$ and non-zero components are {\bf sign-alternating}, that is, the non-zero components follow each other with alternating sign.
  $\bullet$
\end{remark}

\begin{remark}
  The notion of elementary prefix-bounded g-polymatroids can be extended in a natural way to the case when $\varphi (u_0)$ and $\gamma (u_0)$ are not necessarily zero.
  This model includes the polytope of such sign-alternating $(0,\pm 1)$-valued vectors whose first and last non-zero components are arbitrary --- and not necessarily $+1$ and $+1$.
  $\bullet$
\end{remark}

\begin{remark}\label{rem:single_row_suffment}
  A natural question is whether one still obtains a g-polymatroid when bounds are imposed not only on the prefixes but also on the {\bf suffixes} (i.e., on the complements of the prefixes).
  In general, the answer is negative.
  For example, let $U \coloneqq {u_1,u_2,u_3}$, and let the upper bound be $+1$ both for the prefix ${u_1,u_2}$ and for the suffix ${u_2,u_3}$.
  The corresponding polyhedron is not a g-polymatroid.
  However, if the lower and upper bounds for the entire ground set coincide, then the bounds on the prefixes and suffixes together define a base-polyhedron.

  An analogous situation occurs when bounds are imposed on segments rather than on prefixes and suffixes: in general, the resulting polyhedron is not a g-polymatroid.
  If, however, the lower and upper bounds are uniform, then they define a weak pair, and consequently we do obtain a g-polymatroid.
  $\bullet$
\end{remark}

\section{The intersection of two g-polymatroids with a box and a plank}\label{sec:intersection}

The intersection theorem of two g-polymatroids is as follows, see Proposition~V.3.4 and Theorem~V.1.4 in~\cite{frank1988generalized}.
\begin{proposition}\label{prop:metszet}
  {\bf (A)} \ Let $Q_1 \coloneqq Q(p^*_1,b^*_1)$ and $Q_2 \coloneqq Q(p^*_2,b^*_2)$ be two g-polymatroids bordered by integer-valued strong pairs $(p^*_1,b^*_1)$ and $(p^*_2,b^*_2)$.
  Their intersection $Q \coloneqq Q_1\cap Q_2$ is an integral polyhedron which is non-empty if and only if
  \[
    p^*_1\leq b^*_2 \ \ \ \hbox{and}\ \ \ p^*_2\leq b^*_1.
  \]

  \noindent {\bf (B)} \ Let $(p_1,b_1)$ and $(p_2,b_2)$ be integer-valued weak pairs.
  Then the linear inequality system
  \begin{equation}\label{eq:box-tdi}
    \bigl\{\max \{p_1(X), p_2(X) \} \leq \widetilde x(X) \leq \min \{b_1(X), b_2(X)\} \hbox{ for every } X\subseteq S \bigr\}
  \end{equation}
  is box-TDI.
  \FBOX
\end{proposition}
It is easy to see that Part (A) of Proposition~\ref{prop:metszet} is equivalent to the discrete separation theorem.
Since lower and upper bounds $f\leq g$ for components and $\alpha \leq \beta$ for the total sum of components can be incorporated in a weak pair, it follows that the linear inequality system obtained from~\eqref{eq:box-tdi} by adding $f\leq x\leq g$ and $\alpha \leq \widetilde x(S)\leq \beta$ is also box-TDI.
It should be noted that Part (B) was formulated earlier for polymatroids by Edmonds~\cite{edmonds1970submodular}.

We note that the set of integral elements of an integral g-polymatroid (or of the intersection of two integral g-polymatroids) is also known as an {\bf M$\nat$-convex} set (respectively, an {\bf M$\nat_2$-convex} set), as introduced by Murota~\cite{Murota03}, which play a fast increasing role in discrete convex optimization.

\begin{proposition}\label{prop:p0b0f0g0}
  Let $p_0$ be a fully supermodular function on a ground set $S_0$ and let $b_0$ be a fully submodular function, and suppose that $p_0(S_0) = b_0(S_0) = 0$.
  Furthermore, let $f_0$ and $g_0$ be lower and upper bound functions on $S_0$ for which $f_0\leq g_0$.
  The intersection $B_0 \coloneqq B(b_0)\cap B'(p_0)\cap T(f_0,g_0)$ is non-empty if and only if
  \begin{equation}\label{eq:Y1Y2}
    p_0(Y_1) + \widetilde f_0(Y_2-Y_1) \leq b_0(Y_2) + \widetilde g_0(Y_1-Y_2)
  \end{equation}
  holds for every $Y_1,Y_2\subseteq S_0$.
\end{proposition}
\begin{proof}
  For an element $z\in B(b_0)\cap B'(p_0)\cap T(f_0,g_0)$, we have
  \[
    p_0(Y_1) + \widetilde f_0(Y_2-Y_1) \leq \widetilde z(Y_1) + \widetilde z(Y_2-Y_1) = \widetilde z(Y_1\cup Y_2)
  \]
  and
  \[
    b_0(Y_2) + \widetilde g_0(Y_1-Y_2) \geq \widetilde z(Y_2) + \widetilde z(Y_1-Y_2) = \widetilde z(Y_1\cup Y_2),
  \]
  from which~\eqref{eq:Y1Y2} follows.

  Now we prove the reverse direction.
  First, apply Proposition~\ref{prop:j17-2.6} for the following weak pair defining the intersection $B'(p_0) \cap T(f_0,+\infty )$:
  \[
    p(X) \coloneqq
    \begin{cases}
      \max\{\widetilde f_0(X),p_0(X)\} & \ \ \hbox{if}\ \ \ |X| = 1,\\
      p_0(X) & \ \ \hbox{if}\ \ \ |X| \neq 1,
    \end{cases}
    \qquad
    b(X) \coloneqq
    \begin{cases}
      0& \ \ \hbox{if}\ \ \ X = S_0,\\
      +\infty & \ \ \hbox{if}\ \ \ X \neq S_0,
    \end{cases}
  \]
  where $X\subseteq S_0$.
  Observe that the condition of non-emptiness is $p_0(Y_1) + \widetilde f_0(S_0-Y_1) \leq 0$ for every $Y_1\subseteq S_0$, which holds by~\eqref{eq:Y1Y2}.
  Thus, the 0-base-polyhedron $B'(p_0) \cap T(f_0,+\infty)$ is non-empty.

  It follows from Proposition~\ref{prop:j17-2.11} that the unique fully supermodular function $p_0^f$ bordering this $0$-base-polyhedron from below is
  \[
    p^f_0(U) = \max\{p_0(Y) + \widetilde f_0(U-Y) : Y\subseteq U\}.
  \]

  We obtain analogously that the 0-base-polyhedron $B(b_0)\cap T(-\infty ,g)$ is also non-empty, and the unique fully submodular function $b_0^g$ bordering this $0$-base-polyhedron from above is
  \[
    b^g_0(U) = \min\{b_0(Y) + \widetilde g_0(U-Y) : Y\subseteq U\}.
  \]

  It follows from the discrete separation theorem that $B_0$ is non-empty if and only if $p_0 ^f\leq b ^g_0$, which is equivalent to requiring that
  \[
    p_0(Y_1) + \widetilde f_0(U-Y_1) \leq b_0(Y_2) + \widetilde g_0(U-Y_2)
  \]
  holds for every $U\subseteq S_0$ and $Y_1, Y_2\subseteq U$.
  But if this is violated for some $U,Y_1,Y_2$, then $f_0\leq g_0$ implies that it is violated when $U$ is replaced with $U' \coloneqq Y_1\cup Y_2$, that is, $p_0(Y_1) + \widetilde f_0(Y_2-Y_1) > b_0(Y_2) + \widetilde g_0(Y_1-Y_2)$, contradicting~\eqref{eq:Y1Y2}.
\end{proof}

In the literature, there is a characterization for the non-emptiness of the intersection of two g-polymatroids described by a weak pair, see~\cite{frank1984generalized}, but the condition is a bit cumbersome.
Fortunately, we will need only a special case when we intersect two g-polymatroids, given by their strong pairs, with a box and a plank.
\begin{theorem}\label{thm:Q1Q2T}
  We are given two g-polymatroids $Q_1 \coloneqq Q(p^*_1,b^*_1)$ and $Q_2 \coloneqq Q(p^*_2,b^*_2)$ on a ground set $S$ where $(p^*_1,b^*_1)$ and $(p^*_2,b^*_2)$ are strong pairs.
  Furthermore, let $f$ and $g$ be lower and upper bounds on $S$ for which $f\leq g$, where $f$ may take $-\infty$ values and $g$ may take $+\infty$ values.
  Let $T \coloneqq T(f,g) = \{x : f\leq x\leq g\}$ denote the box defined by $f$ and $g$.
  Finally, let $\alpha$ and $\beta$ be numbers serving as lower and upper bounds for which $-\infty \leq \alpha \leq \beta \leq +\infty$, and let $K = K(\alpha,\beta )$ denote the plank defined by $\alpha$ and $\beta$.

  Then the intersection $Q_1\cap Q_2\cap T\cap K$ is non-empty if and only if each of
  \begin{gather}
    p^*_1(X_1) + \widetilde f(X_2-X_1) \leq b^*_2(X_2) + \widetilde g(X_1-X_2),\label{eq:Q1Q2Ta}\\
    p^*_2(X_2) + \widetilde f(X_1-X_2) \leq b^*_1(X_1) + \widetilde g(X_2-X_1),\label{eq:Q1Q2Tb}\\
    b^*_1(X_1) + b^*_2(X_2) + \widetilde g(\ol{X_1} \cap \ol{X_2}) - \widetilde f(X_1\cap X_2) \geq \alpha ,\label{eq:Q1Q2Talfa}\\
    p^*_1(X_1) + p^*_2(X_2) + \widetilde f(\ol{X_1} \cap \ol{X_2}) - \widetilde g(X_1\cap X_2) \leq \beta\label{eq:Q1Q2Tbeta}
  \end{gather}
  holds for every $X_1,X_2\subseteq S$.
\end{theorem}
\begin{proof}
  To prove necessity, let $z$ be an element of the intersection $Q_1\cap Q_2\cap T\cap K$.
  Then
  \[
    p^*_1(X_1) +\widetilde f(X_2-X_1) \leq \widetilde z_1(X_1) +\widetilde z(X_2-X_1) = \widetilde z(X_2) +\widetilde z(X_1-X_2) \leq b^*_2(X_2) +\widetilde g(X_1-X_2),
  \]
  from which~\eqref{eq:Q1Q2Ta} follows.
  The necessity of~\eqref{eq:Q1Q2Tb} follows similarly.
  Moreover, it is true for $z$ that
  \begin{align*}
    \alpha &\leq \widetilde z(S) = \widetilde z(X_1) + \widetilde z(X_2) -\widetilde z(X_1\cap X_2) + \widetilde z(\ol{X_1} \cap \ol{X_2})\\
           &\leq b^*_1(X_1) + b^*_2(X_2) - \widetilde f(X_1\cap X_2) + \widetilde g(\ol{X_1} \cap \ol{X_2}),
  \end{align*}
  from which~\eqref{eq:Q1Q2Talfa} follows.
  The necessity of~\eqref{eq:Q1Q2Tbeta} can be seen analogously.
  \medskip

  To prove sufficiency, let $s_0$ be a new element (not occurring in $S$) and let $S_0 \coloneqq S+s_0$.
  Define the bounding functions $f_0$ and $g_0$ on $S_0$ as
  \[
    f_0(s) \coloneqq
    \begin{cases}
      f(s) & \ \ \hbox{if}\ \ \ s\in S,\\
      -\beta & \ \ \hbox{if}\ \ \ s = s_0,
    \end{cases}
    \qquad
    g_0(s) \coloneqq
    \begin{cases}
      g(s) & \ \ \hbox{if}\ \ \ s\in S,\\
      -\alpha & \ \ \hbox{if}\ \ \ s = s_0,
    \end{cases}
  \]
  and let $T_0 \coloneqq T(f_0,g_0).$
  Define two set functions $p^*_0$ and $b^*_0$ on $S_0$ as
  \[
    p^*_0(Y) \coloneqq
    \begin{cases}
      p^*_1(Y) & \ \ \hbox{if}\ \ \ s_0 \notin Y,\\
      -b^*_1(S_0-Y) & \ \ \hbox{if}\ \ \ s_0\in Y,
    \end{cases}
    \qquad
    b^*_0(Y) \coloneqq
    \begin{cases}
      b^*_2(Y) & \ \ \hbox{if}\ \ \ s_0 \notin Y,\\
      -p^*_2(S_0-Y) & \ \ \hbox{if}\ \ \ s_0\in Y.
    \end{cases}
  \]
  Now $p^*_0$ is fully supermodular for which $p^*_0(S_0) = 0$, and $b^*_0$ is fully submodular for which $b^*_0(S_0) = 0$.
  It follows readily that for any element $z \in Q_1 \cap Q_2 \cap T \cap K$, the vector $z_0$ defined as
  \[
    z_0(s) \coloneqq
    \begin{cases}
      z(s) & \ \ \hbox{if}\ \ \ s\in S,\\
      -\widetilde z(S) & \ \ \hbox{if}\ \ \ s = s_0
    \end{cases}
  \]
  lies in $B'(p_0^*)\cap B(b_0^*) \cap T_0$.
  Conversely, if $z_0 \in B'(p_0^*)\cap B(b_0^*)\cap T_0$, then the restriction of $z_0$ to $S$ lies in $Q_1\cap Q_2\cap T\cap K$.
  That is, the original intersection $Q_1\cap Q_2\cap T\cap K$ is non-empty if and only if the intersection $B'(p^*_0) \cap B(b^*_0) \cap T_0$ is non-empty.

  Apply Proposition~\ref{prop:p0b0f0g0} to $p^*_0$ and $b^*_0$ in place of $p_0$ and $b_0$, and assume indirectly that~\eqref{eq:Y1Y2}~is violated by $Y_1$ and $Y_2$, that is,
  \begin{equation}\label{eq:kontra}
    p^*_0(Y_1) + \widetilde f_0(Y_2-Y_1)> b^*_0(Y_2) + \widetilde g_0(Y_1-Y_2).
  \end{equation}

  If $Y_1,Y_2\subseteq S$, then, for $X_1 \coloneqq Y_1$ and $X_2 \coloneqq Y_2$,~\eqref{eq:kontra} means

  \[
    p^*_1(X_1) + \widetilde f(X_2-X_1) > b^*_2(X_2) + \widetilde g(X_1-X_2),
  \]
  contradicting~\eqref{eq:Q1Q2Ta}.

  If $s_0\in Y_1\cap Y_2$, then, for $X_1 \coloneqq S_0-Y_1$ and $X_2 \coloneqq S_0-Y_2$,~\eqref{eq:kontra} means
  \[
    p^*_2(X_2) + \widetilde f(X_1-X_2) > b^*_1(X_1) + \widetilde g(X_2-X_1),
  \]
  contradicting~\eqref{eq:Q1Q2Tb}.

  If $s_0\in Y_1$ and $Y_2\subseteq S$, then~\eqref{eq:kontra} implies
  \[
    -b^*_1(S_0-Y_1) + \widetilde f_0(Y_2-Y_1) > b^*_2(Y_2) + \widetilde g_0(Y_1-Y_2).
  \]
  In this case, for $X_1 \coloneqq S_0-Y_1$ and $X_2 \coloneqq Y_2$, we have $Y_2-Y_1 = X_1\cap X_2$, \ $Y_1-Y_2 = S-(X_1\cup X_2) + s_0 = \ol{X_1} \cap \ol{X_2} + s_0$, moreover, $g(s_0) = -\alpha$, and hence
  \[
    -b^*_1(X_1) + \widetilde f(X_1\cap X_2) > b^*_2(X_2) + \widetilde g(\ol{X_1} \cap \ol{X_2}) - \alpha,
  \]
  contradicting~\eqref{eq:Q1Q2Talfa}.

  Finally, if $Y_1\subseteq S$ and $s_0\in Y_2$, then it follows from~\eqref{eq:kontra} that
  \[
    p^*_1(Y_1) + \widetilde f_0(Y_2-Y_1) > - p^*_2(S_0-Y_2) + \widetilde g_0(Y_1-Y_2).
  \]
  In this case, for $X_1 \coloneqq Y_1$ and $X_2 \coloneqq S_0- Y_2$, we have $Y_1-Y_2 =  X_1\cap X_2$, \ $Y_2-Y_1 = S- (X_1\cup X_2) + s_0 = \ol{X_1} \cap \ol{X_2} +s_0$, moreover, $f(s_0) = -\beta$, and hence
  \[
    p^*_1(X_1) + \widetilde f(\ol{X_1} \cap \ol{X_2}) - \beta > -p^*_2(X_2) + \widetilde g(X_1\cap X_2),
  \]
  contradicting~\eqref{eq:Q1Q2Tbeta}.
\end{proof}

\begin{corollary}\label{cor:Bresz}
  If both g-polymatroids in Theorem~\ref{thm:Q1Q2T} are base-polyhedra, that is, $Q_1 \coloneqq B(b^*_1)$ and $Q_2 \coloneqq B(b^*_2)$ for which $H \coloneqq b^*_1(S) = b^*_2(S)$, and $\alpha = -\infty,\ \beta = +\infty$, then~\eqref{eq:Q1Q2Ta} and~\eqref{eq:Q1Q2Tb} are equivalent, and hence, the intersection $Q_1\cap Q_2\cap T\cap K$ is non-empty if and only if
  \begin{equation}\label{eq:ekvi}
    H \leq b^*_1(X') + b^*_2(X'') + \widetilde g(\ol{X'} \cap \ol{X''}) - \widetilde f(X' \cap X'')
  \end{equation}
  holds for every $X',X''\subseteq S$.
\end{corollary}
\begin{proof}
  Apply Theorem~\ref{thm:Q1Q2T} to the case when $\alpha = -\infty,\ \beta = +\infty$.
  In this case, only~\eqref{eq:Q1Q2Ta} and~\eqref{eq:Q1Q2Tb} matter.
  Observe that for $X' \coloneqq \ol{X_1}$ and $X'' \coloneqq X_2$, we have $ X_2-X_1 = X'\cap X''$, $X_1-X_2 = \ol{X'} \cap \ol{X''}$, $b^*_2(X_2) = b^*_2(X'')$, and Claim~\ref{cl:pSbS} implies that $p^*_1(X_1) = H- b^*_1(X')$.
  Therefore, the inequalities in~\eqref{eq:Q1Q2Ta} and~\eqref{eq:ekvi} are equivalent.
  We get in an analogous way that, for $X' \coloneqq X_1$ and $X'' \coloneqq \ol{X_2}$, the inequalities in~\eqref{eq:Q1Q2Tb} and~\eqref{eq:ekvi} are equivalent.
\end{proof}

The following corollary can immediately be obtained from Theorem~\ref{thm:Q1Q2T}.
\begin{corollary}
  Let $Q_1 \coloneqq Q(p^*_1,b^*_1), Q_2 \coloneqq  Q(p^*_2,b^*_2)$, and $T(f,g)$ be as in Theorem~\ref{thm:Q1Q2T}, and suppose that $Q \coloneqq Q_1 \cap Q_2\cap T(f,g)$ is non-empty.
  Then
  \begin{gather*}
    \max \{\widetilde x(S) : x\in Q\} = \min \{b^*_1(X_1) + b^*_2(X_2) + \widetilde g(\ol{X_1} \cap \ol{X_2}) - \widetilde f(X_1\cap X_2) : X_1, X_2\subseteq S\},\\
    \min \{\widetilde x(S) : x\in Q\} = \max \{ p^*_1(X_1) + p^*_2(X_2) + \widetilde f(\ol{X_1} \cap \ol{X_2}) - \widetilde g(X_1\cap X_2) : X_1, X_2\subseteq S\}.
  \end{gather*}
  \FBOX
\end{corollary}

To close this section, we investigate the special case when a specified finite-valued prescription vector $m_P$ is defined on a subset $S_P$ of $S$, and we are interested in the existence of an element in the intersection of two g-polymatroids whose components in $S_P$ are the prescribed values.
\begin{corollary}\label{cor:fixed}
  Let $Q_1 \coloneqq Q(p^*_1,b^*_1)$ and $Q_2 \coloneqq Q(p^*_2,b^*_2)$ be g-polymatroids in $\R^S$ with strong bordering pairs $(p^*_1,b^*_1)$ and $(p^*_2,b^*_2)$.
  Furthermore, let $m_P$ be a finite-valued vector on a subset $S_P$ of $S$.

  {\bf (A)} \ The intersection $Q_1\cap Q_2$ has an integral element $z$ for which $z\vert S_P = m_P$ (that is, the restriction of $z$ to $S_P$ equals $m_P$) if and only if the following two inequalities hold for all subsets $X_1, X_2 \subseteq S$ with $X_1 \ominus X_2 \subseteq S_P$:
  \begin{gather}
    p^*_1(X_1) + \widetilde m_P(X_2-X_1) \leq b^*_2(X_2) + \widetilde m_P(X_1-X_2),\label{eq:Q1Qq-mF-1}\\
    p^*_2(X_2) + \widetilde m_P(X_1-X_2) \leq b^*_1(X_1) + \widetilde m_P(X_2-X_1).\label{eq:q1qq-mf-2}
  \end{gather}

  {\bf (B)} \ In the special case when $Q_1$ and $Q_2$ are base-polyhedra for which $H \coloneqq b^*_1(S) = b^*_2(S)$, there exists an integral element $z\in Q_1 \cap Q_2$ for which $z\vert S_P = m_P$ if and only if
  \begin{equation}\label{eq:ekvib}
    H \leq b^*_1(X') + b^*_2(X'') + \widetilde m_P(\ol{X'} \cap \ol{X''}) - \widetilde m_P(X' \cap X'')
  \end{equation}
  holds for every $X',X'' \subseteq S$ for which $\ol{S_P} \subseteq X'\ominus X''$.
\end{corollary}
\begin{proof}
  Let $\alpha \coloneqq -\infty,\ \beta \coloneqq +\infty$, and let
  \[
    f(s) \coloneqq
    \begin{cases}
      m_P(s) & \ \ \hbox{if}\ \ \ s\in S_P,\\
      -\infty & \ \ \hbox{if}\ \ \ s\in \ol{S_P},
    \end{cases}
    \qquad
    g(s) \coloneqq
    \begin{cases}
      m_P(s) & \ \ \hbox{if}\ \ \ s\in S_P,\\
      +\infty & \ \ \hbox{if}\ \ \ s\in \ol{S_P}.
    \end{cases}
  \]

  To prove Part (A), observe that condition~\eqref{eq:Q1Q2Ta} holds automatically when $\widetilde f(X_2-X_1) = -\infty$ or $\widetilde g(X_1-X_2) = +\infty$, therefore, it suffices to require~\eqref{eq:Q1Q2Ta} only when $X_1 \ominus X_2 \subseteq S_P$, which is just~\eqref{eq:Q1Qq-mF-1}.
  Analogously, in the present special case,~\eqref{eq:Q1Q2Tb} is~\eqref{eq:q1qq-mf-2}.

  To prove Part B, observe that~\eqref{eq:ekvi} holds automatically when $\widetilde f(X'\cap X'') = -\infty$ or ${\widetilde g(\ol{X'} \cap \ol{X''}) = +\infty}$, that is, it suffices to require~\eqref{eq:ekvi} only when $X'\cap X''\subseteq S_P$ and $\ol{X'} \cap \ol{X''}\subseteq S_P$, that is, $\ol{S_P}\subseteq X'\ominus X''$.
  In this case,~\eqref{eq:ekvi} and~\eqref{eq:ekvib} are equivalent.
\end{proof}

\begin{remark}
  The concept of submodular flows was first introduced and investigated by Edmonds and Giles~\cite{edmonds1977min}.
  It is a known fact that the intersection of two g-polymatroids with a box and a plank can be obtained as the projection of a submodular flow polyhedron, and such a projection is also known to be a submodular flow polyhedron.

  Since there are (combinatorial) strongly polynomial algorithms for deciding whether a submodular flow polyhedron is empty or not, and, in the latter case, also for computing a cheapest submodular flow, it follows that the problems concerning the intersection of two g-polymatroids are algorithmically tractable.
  In the next sections, we consider only the special case when both g-polymatroids are laminar.
  In this case, the feasibility and optimization problems can be solved via standard network-flow algorithms, see Section~\ref{sec:circulation}.
  $\bullet$
\end{remark}

It is known that while a single g-polymatroid admits the integer decomposition property, the same does not generally hold for the intersection of two g-polymatroids.
However, in the special case of laminar g-polymatroids, the situation is more favorable.
Recall that the intersection of a laminar g-polymatroid with a box and a plank yields another laminar g-polymatroid, and that the intersection of two laminar g-polymatroids can be described by a linear inequality system with a totally unimodular constraint matrix.
Therefore, Proposition~\ref{prop:strong-Cara} directly yields the following result.
\begin{corollary}\label{cor:integer-Cara}
  The intersection of two laminar g-polymatroids (in particular, prefix-bounded g-polymatroids) with a box and a plank admits the sign-consistent integer Carath\'eodory property, and hence also the sign-consistent integer decomposition property.
  \FBOX
\end{corollary}

Note that a sign-consistent integer decomposition can be found in strongly polynomial time via standard network-flow algorithms. 


\section{Prefix-bounded matrices}\label{sec:PBM}

Heuer and Striker~\cite{heuer2022partial} introduced the notion of a {\bf partial alternating sign matrix} (PASM), defined as a $(0, \pm1)$-valued $m \times n$ matrix in which the sum of entries in every horizontal and vertical prefix (i.e., prefix of a row and column, respectively) is either $0$ or $+1$.
Throughout, we use the term {\bf prefix} to refer to both horizontal and vertical prefixes, specifying the type only when necessary for clarity.
It is straightforward to verify that any upper-left rectangular submatrix of an ASM is a PASM, and conversely, every PASM can be obtained this way.
Heuer and Striker~\cite{heuer2022partial} showed that the convex hull of $m \times n$ PASMs can be described as $\{ x \in \R^{m \times n} : 0 \leq \widetilde{x}(Z) \leq 1 \text{ for every prefix } Z\}$, where $\widetilde{x}(Z)$ denotes the sum of entries in prefix $Z$, see Theorem~4.6 of their paper.

Another extension of ASMs was investigated by Behrend and Knight~\cite{behrend2007higher}.
They call an integer-valued square matrix a {\bf higher spin alternating sign matrix} if, for a positive integer $r$, the sum of entries of all prefixes is non-negative and at most $r$, and this sum is exactly $r$ for all full rows and columns.
For $r = +1$, we get back the class of ASMs.
For larger values of $r$, however, unlike classical ASMs, the signs of non-zero entries in a row or column of a higher spin ASM do not necessarily alternate.
Note that the polytope of higher spin alternating sign matrices is precisely the $r$-elongation of the polytope of ASMs.

The notion of higher spin alternating matrices was further generalized by Brualdi and Dahl~\cite{brualdi2017alternatingExtensions}.
Given an $n \times n$ matrix $B \coloneqq (b_{i,j})$ of non-negative integers, an integer-valued $n \times n$ matrix is said to be {\bf sum-majorized} by $B$ if, for each $(i,j)$, both the sum of entries in the horizontal prefix ending at $(i,j)$ and the sum of entries in the vertical prefix ending at $(i,j)$ are non-negative and at most $b_{i,j}$.
Additionally, the sum of entries in row $i$ must be exactly $b_{i,n}$ for all $i = 1, \dots, n$, and the sum of entries in column $j$ must be exactly $b_{n,j}$ for all $j = 1, \dots, n$.

A related development appears in a more recent paper by Brualdi and Dahl~\cite{brualdi2023multi}, where they study $(0,\pm 1)$-valued $m \times n$ matrices such that, for given non-negative integers $r_1, \dots, r_m$ and $s_1, \dots, s_n$, the sum of entries in row $i$ is required to be exactly $r_i$, and for every proper prefix of that row, the sum must lie between $0$ and $r_i$, and analogously, the sum of entries in column $j$ is exactly $s_j$, and for every proper prefix of that column, the sum is constrained to lie between $0$ and $s_j$.
Among other results, the authors provide a polyhedral description of the convex hull of these matrices.

Yet another extension of ASMs was introduced by Aval~\cite{aval2008keys}, who defined a $(0, \pm1)$-valued $m \times n$ matrix to be a {\bf sign matrix} if the sum of entries in every vertical prefix is either $0$ or $+1$, and the sum of entries in every horizontal prefix is non-negative.
For further studies of related classes of sign-alternating matrices, see the work of Heuer, Solhjem, and Striker~\cite{heuer2024nu, heuer2022partial, solhjem2019sign}.

All of the matrix classes discussed above can be viewed as special cases of the following general framework.
We call an integer-valued $m \times n$ matrix a {\bf prefix-bounded matrix} (PBM) if the sum of entries in each horizontal and vertical prefix lies between specified lower and upper bounds, which may vary from prefix to prefix.
Our goal is to establish a necessary and sufficient condition for the existence of such a matrix, subject also to given lower and upper bounds on the individual entries, and on the total sum of the entries.
Figure~\ref{fig:pbm} illustrates these constraints.
\begin{figure}[h]
  \centering
  \begin{tikzpicture}[scale=1.5]
    \draw[step=1cm,black!80,thin] (0,0) grid (3,3);

    \foreach \i in {0,...,2} {
      \draw[blue!90,opacity=0.45,line width=28pt,line cap=round]
      (0.5,\i + 0.5) -- (2.5,\i + 0.5);
    }
    \foreach \j in {0,...,2} {
      \draw[red!90,opacity=0.45,line width=28pt,line cap=round]
      (\j + 0.5, 2.5) -- (\j + 0.5, 0.5);
    }
    \foreach \i in {0,...,2} {
      \draw[blue!90,opacity=0.6,line width=18pt,line cap=round]
      (0.5,\i + 0.5) -- (1.5,\i + 0.5);
    }
    \foreach \j in {0,...,2} {
      \draw[red!90,opacity=0.6,line width=18pt,line cap=round]
      (\j + 0.5, 2.5) -- (\j + 0.5, 1.5);
    }
    \foreach \x in {1.5,2.5} {
      \fill[red] (\x,2.5) circle (3pt);
    }
    \foreach \y in {1.5,0.5} {
      \fill[blue] (0.5,\y) circle (3pt);
    }
    \begin{scope}
      \clip (0.5,2.5) circle (3pt); 
      \fill[blue!90] (0.5,2.5) -- ++(135:0.2) arc (135:315:0.2) -- cycle;
    \end{scope}
    \begin{scope}
      \clip (0.5,2.5) circle (3pt);
      \fill[red!90] (0.5,2.5) -- ++(315:0.2) arc (315:495:0.2) -- cycle;
    \end{scope}

  \end{tikzpicture}
  \caption{Visualization of a $3 \times 3$ grid with all row prefixes (blue horizontal bars) and column prefixes (red vertical bars).
    In a prefix-bounded matrix, the sum in each such prefix is subject to lower and upper bounds.
    Additional constraints may also apply to individual entries and to the total sum of all entries.}\label{fig:pbm}
\end{figure}

Recall that $S_{m,n} \coloneqq [m] \times [n]$ denotes the set of positions of the entries in an $m \times n$ matrix.
We are given two integer-valued functions $f\leq g$ on $S_{m,n}$ serving as lower and upper bounds for the entries.
Moreover, we are given a lower bound $\alpha$ and an upper bound $\beta$ for the total sum of the entries, we assume that $\alpha \leq \beta$.
Here $f$ and $\alpha$ may take $-\infty$ values, $g$ and $\beta$ may take $+\infty$ values.
We are also given four integer-valued $m \times n$ matrices $\Phi^1$, $\Phi^2$, $\Gamma^1$, and $\Gamma^2$, where $\Phi^1$ and $\Phi^2$ may have entries of value $-\infty$ and $\Gamma^1$ and $\Gamma^2$ may have entries of value $+\infty$.
These matrices specify lower and upper bounds on the entry sums of horizontal and vertical prefixes, respectively.
For example $\Phi^1(i,j)$ is the lower bound for the sum of the entries in the horizontal prefix ending at position $(i,j)$.
We assume that $\Phi^1 \leq \Gamma^1$ and $\Phi^2 \leq \Gamma^2$.

Recall from Section~\ref{sec:pb-g-poli} that given a totally ordered set $U = \{u_1, \dots, u_t\}$ and functions $\varphi$ and $\gamma$ on $U$, we defined the set functions $p^*$ and $b^*$ as in equations~\eqref{eq:p'b'def} and~\eqref{eq:addit}.
We now apply this separately to each row and each column of $S_{m,n}$ with the natural left-to-right or top-to-bottom ordering as the total order, respectively.

First, consider the rows.
For the $i$-th row ($i = 1,\dots ,m$), let $Q^1_i$ denote the prefix-bounded (and thus laminar) g-polymatroid defined by the $i$-th row of the matrices $\Phi^1$ and $\Gamma^1$.
In this context, the underlying ground set $U = \{u_1,\dots ,u_t\}$ corresponds to the positions of the $i$-th row, hence $t = n$, with $\varphi(u_k) = \Phi^1 ({i,k})$ and $\gamma(u_k) = \Gamma^1 ({i,k})$ for $k = 1, \dots, n$.
Let $Q^1$ be the direct sum of these $m$ g-polymatroids.
Let $P^1_{i,j}$ denote the horizontal prefix of $S_{m,n}$ consisting of the first $j$ positions of the $i$-th row.
Then
\begin{equation}\label{eq:Q1def}
  Q^1 = \{x\in \R^{m \times n} : \Phi^1 (i,j) \leq \widetilde x(P^1_{i,j}) \leq \Gamma^1 (i,j) \hbox{ for every } (i,j) \in S_{m,n} \}.
\end{equation}

Second, we define the prefix-bounded g-polymatroid $Q^2$ analogously as the direct sum of the $n$ elementary prefix-bounded g-polymatroids on the columns of $S_{m,n}$.
Let $P^2_{i,j}$ denote the vertical prefix of $S_{m,n}$ consisting of the first $i$ positions of the $j$-th column.
Then
\begin{equation}\label{eq:Q2def}
  Q^2 = \{x\in \R^{m \times n} : \Phi^2 (i,j) \leq \widetilde x(P^2_{i,j}) \leq \Gamma^2 (i,j) \hbox{ for every } (i,j) \in S_{m,n} \}.
\end{equation}

The unique strong pairs $(p^*_1,b^*_1)$ and $(p^*_2, b^*_2)$ bordering $Q^1$ and $Q ^2$ can be obtained by applying Lemma~\ref{lem:p*b*} separately to each row and each column, respectively, which we will derive later in this section.

With these constructions in place, and invoking Proposition~\ref{prop:metszet}, we arrive at the following result.
\begin{proposition}
  The convex hull of prefix-bounded matrices defined by the horizontal-prefix bounds $(\Phi^1, \Gamma^1)$ and vertical-prefix bounds $(\Phi^2, \Gamma^2)$ is the intersection of the prefix-bounded g-polymatroids $Q^1$ and $Q^2$ given in~\eqref{eq:Q1def} and~\eqref{eq:Q2def}.
  In particular, the set of these prefix-bounded matrices is an M$\nat _2$-convex set.
  \FBOX
\end{proposition}

From this and from Theorem~\ref{thm:box-tdi1}, we obtain the following two corollaries.
\begin{corollary}\label{cor:box-TDI}
  The linear inequality system
  \[
    \begin{cases}
      & \Phi^1 (i,j) \leq \widetilde x( P^1_{i,j}) \leq \Gamma^1 (i,j) \hbox{ for every } (i,j) \in S_{m,n}\\
      & \Phi^2 (i,j) \leq \widetilde x(P^2_{i,j}) \leq \Gamma^2 (i,j) \hbox{ for every } (i,j) \in S_{m,n}
    \end{cases}
  \]
  describing the convex hull of prefix-bounded matrices defined by $\Phi^1, \Gamma^1, \Phi^2,$ and $\Gamma^2$ is box-TDI.
  \FBOX
\end{corollary}

\begin{corollary}\label{cor:boundedPBMPolytope}
  For given prefix bounds $\Phi^1 \leq \Gamma^1$ and $\Phi^2 \leq \Gamma^2$, the convex hull of the prefix-bounded matrices of size $m \times n$ whose entries lie between $f$ and $g$, and whose total sum lies between $\alpha$ and $\beta$ is described by the linear inequality system
  \[
    \begin{cases}
      & \Phi^1 (i,j) \leq \widetilde x( P^1_{i,j}) \leq \Gamma^1 (i,j) \hbox{ for every } (i,j) \in S_{m,n}\\
      & \Phi^2 (i,j) \leq \widetilde x(P^2_{i,j}) \leq \Gamma^2 (i,j) \hbox{ for every } (i,j) \in S_{m,n}\\
      & f \leq x \leq g\\
      & \alpha \leq \widetilde x(S_{m,n}) \leq \beta.
    \end{cases}
  \]
  \FBOX
\end{corollary}

This polyhedron is the intersection of two laminar g-polymatroids, as it is the intersection of the (prefix-bounded) g-polymatroids $Q^1$ and $Q^2$ given in~\eqref{eq:Q1def} and~\eqref{eq:Q2def} with the box $T(f,g)$ and the plank $K(\alpha, \beta)$.
As we noted in Section~\ref{sec:lam.g-poli}, the linear optimization problem concerning the intersection of two laminar g-polymatroids can be formulated as an optimal feasible circulation problem.
In the special case of prefix-bounded matrices, the details of this approach will be worked out in Section~\ref{sec:circulation}.

As a special case of Corollary~\ref{cor:integer-Cara}, we get the following.
\begin{corollary}\label{cor:int-decomp-spec}
  The convex hull of prefix-bounded matrices with bounds on their entries and total sum admits the sign-consistent integer Carath\'eodory property --- and therefore also the sign-consistent integer decomposition property.
  In particular, for given prefix bounds $\Phi^1 \leq \Gamma^1$ and $\Phi^2 \leq \Gamma^2$, let $A$ be an $m \times n$ prefix-bounded matrix whose entries are between $f$ and $g$, and for which the total sum of entries lies between $\alpha$ and $\beta$.
  For every positive integer $k$, there exist $m \times n$ prefix-bounded matrices $A_1, \dots, A_k$ defined by $\left\lfloor\frac{\Phi^1}{k}\right\rfloor$, $\Bigl\lceil\frac{\Gamma^1}{k}\Bigr\rceil$, $\left\lfloor\frac{\Phi^2}{k}\right\rfloor$, and $\Bigl\lceil\frac{\Gamma^2}{k}\Bigr\rceil$ whose entries are between $\left\lfloor\frac{f}{k}\right\rfloor$ and $\Bigl\lceil\frac{g}{k}\Bigr\rceil$, and for which the total sum of entries lies between $\Bigl\lfloor\frac{\alpha}{k}\Bigr\rfloor$ and $\left\lceil\frac{\beta}{k}\right\rceil$, such that $A = A_1 + \dots + A_k$ and $A_1, \dots, A_k$ are sign-consistent with $A$.
  \FBOX
\end{corollary}

Note that the integer Carath\'eodory property implies that the family $\mathcal{A} = \{A_1, \dots, A_k\}$ in Corollary~\ref{cor:int-decomp-spec} can be chosen such that the number of distinct members of the family $\mathcal{A}$ is at most $nm+1$.

Corollary~\ref{cor:int-decomp-spec} immediately yields a result conjectured by Brualdi and Dahl~\cite{brualdi2023multi}, formulated as Conjecture~5.12 in their paper, which we now describe.
For a positive integer $k$, they call a $(0,\pm1)$-valued square matrix a {\bf $\boldsymbol{k}$-regular ASM} if the entry sum of each row and each column is $k$, and the entry sum of each horizontal and vertical prefix lies between $0$ and $k$.
\begin{corollary}\label{cor:B+D-conj}
  Every $k$-regular alternating sign matrix can be expressed as the sum of $k$ pairwise pattern-disjoint alternating sign matrices, where the pattern of a matrix is the set of positions of its non-zero entries.
\end{corollary}
\begin{proof}
  The polytope of $k$-regular ASMs is the polyhedron of entry-bounded PBMs in the special case when the lower bound $\Phi^1$ is identically $k$ in the last column and $0$ elsewhere; $\Phi^2$ is identically $k$ in the last row and $0$ elsewhere; the upper bounds $\Gamma^1$ and $\Gamma^2$ are identically $k$; and the entry bounds $f$ and $g$ are identically $-1$ and $+1$, respectively.
  Apply Corollary~\ref{cor:int-decomp-spec} to this setting, and observe that sign-consistency ensures that the obtained ASMs are pairwise pattern-disjoint, completing the proof.
\end{proof}

Now, we turn to the problem of characterizing the existence of PBMs with entry bounds and total-sum bounds.
We first derive the unique strong pairs $(p^*_1, b^*_1)$ and $(p^*_2, b^*_2)$ bordering $Q^1$ and $Q^2$, respectively.
Recall the definition of separated segments within a totally ordered set $U$ from Section~\ref{sec:pb-g-poli}, which we now apply to the rows and columns of the ground set $S_{m,n}$.
There is a natural one-to-one correspondence between subsets of $S_{m,n}$ and systems of separated horizontal segments~---~where no restrictions are imposed between segments of different rows.
For a given subset $X \subseteq S_{m,n}$, let $\mathcal{I}^{\mathrm{ho}}(X)$ and $\mathcal{I}^{\mathrm{ve}}(X)$ denote the system of the maximal horizontal and vertical segments in $X$, respectively.
Both systems are, by definition, separated.

By applying Lemma~\ref{lem:p*b*} separately for each row and column of $S_{m,n}$, we get that the strong pairs $(p^*_1, b^*_1)$ and $(p^*_2, b^*_2)$ bordering $Q^1$ and $Q^2$ are
\begin{equation}\label{eq:bi*pi*}
  \begin{aligned}
    &\begin{cases}
      & p^*_1 (X) 
        = \sum\ [ \Phi^1 (i, j_2) - \Gamma^1 (i, j_1-1) :  \{(i,j_1),\dots,(i,j_2)\} \in \mathcal{I}^{\mathrm{ho}}(X)],\\
      & b^*_1 (X) 
        = \sum\ [ \Gamma^1 (i, j_2) - \Phi^1 (i, j_1-1) :  \{(i,j_1),\dots,(i,j_2)\} \in \mathcal{I}^{\mathrm{ho}}(X)],
    \end{cases}\\
    &\begin{cases}
      & p^*_2 (X) 
        = \sum\ [ \Phi^2 (i_2, j) - \Gamma^2 (i_1-1, j) :  \{(i_1,j),\dots,(i_2,j)\} \in \mathcal{I}^{\mathrm{ve}}(X)],\\
      & b^*_2 (X) 
        = \sum\ [ \Gamma^2 (i_2, j) - \Phi^2 (i_1-1, j) :  \{(i_1,j),\dots,(i_2,j)\} \in \mathcal{I}^{\mathrm{ve}}(X)].
    \end{cases}
  \end{aligned}
\end{equation}

Note that, for any subset $X \subseteq S_{m,n}$, the functions $p ^*_1(X)$ and $b^*_1(X)$ are the strictest lower and upper bounds, respectively, implied by the horizontal-prefix bounds $\Phi^1$ and $\Gamma^1$; analogously, the functions $p ^*_2(X)$ and $b^*_2(X)$ are the strictest lower and upper bounds, respectively, implied by the vertical-prefix bounds $\Phi ^2$ and $\Gamma ^2$.
In what follows, we investigate the existence of prefix-bounded matrices (PBMs), even under additional restrictions.
The functions $p^*_1$, $b^*_1$, $p^*_2$, and $b^*_2$ will be useful for stating the necessary and sufficient conditions.

Theorem~\ref{thm:Q1Q2T} yields the following characterization for the existence of PBMs with entry bounds and total-sum bounds.
\begin{theorem}\label{thm:general1}
  For given prefix bounds $\Phi^1 \leq \Gamma^1$ and $\Phi^2 \leq \Gamma^2$, there exists a prefix-bounded matrix of size $m \times n$ whose entries lie between $f$ and $g$, and for which the total sum of the entries lies between $\alpha$ and $\beta$, if and only if the following inequalities hold for every $X_1, X_2 \subseteq S_{m,n}$:
  \begin{gather}
    p^*_1(X_1) + \widetilde f (X_2-X_1) \leq b^*_2(X_2) + \widetilde g(X_1-X_2),\label{eq:gen1a}
    \\ p^*_2(X_2) + \widetilde f (X_1-X_2) \leq b^*_1(X_1) + \widetilde g(X_2-X_1),\label{eq:gen1b}
    \\ b^*_1(X_1) + b^*_2(X_2) + \widetilde g(\ol{X_1} \cap \ol{X_2}) - \widetilde f(X_1\cap X_2) \geq \alpha,\label{eq:gen1alfa}
    \\ p^*_1(X_1) + p^*_2(X_2) + \widetilde f(\ol{X_1} \cap \ol{X_2}) - \widetilde g(X_1\cap X_2) \leq \beta,\label{eq:gen1beta}
  \end{gather}
  where $p^*_1$, $b^*_1$, $p^*_2$, and $b^*_2$ are the functions defined in~\eqref{eq:bi*pi*}.
  \FBOX
\end{theorem}

Section~\ref{sec:circulation} provides an alternative, self-contained algorithmic proof of this theorem using network-flow techniques instead of g-polymatroids.
We now turn to some useful corollaries of Theorem~\ref{thm:general1}.

\begin{corollary}\label{cor:maxpkm}
  Suppose that, for given prefix bounds $\Phi^1 \leq \Gamma^1$ and $\Phi^2 \leq \Gamma^2$, there exists a prefix-bounded matrix whose entries are between $f$ and $g$.
  Then the maximum of the total sum of entries over such matrices is equal to the value
  \[
    \min \{b^*_1(X_1) + b^*_2(X_2) + \widetilde g(\ol{X_1} \cap \ol{X_2}) - \widetilde f(X_1\cap X_2) : X_1, X_2\subseteq S_{m,n}\},
  \]
  while the minimum of the total sum of entries is equal to
  \[
    \max \{ p^*_1(X_1) + p^*_2(X_2) + \widetilde f(\ol{X_1} \cap \ol{X_2}) - \widetilde g(X_1\cap X_2) : X_1, X_2\subseteq S_{m,n}\},
  \]
  where $p^*_1$, $b^*_1$, $p^*_2$, and $b^*_2$ are the functions defined in~\eqref{eq:bi*pi*}.
  \FBOX
\end{corollary}

In the special case $f\equiv -\infty$, $g\equiv +\infty$, $\alpha = -\infty$, $\beta = +\infty$, Theorem~\ref{thm:general1} reduces to the following form.
\begin{corollary}\label{cor:nofg}
  There exists a prefix-bounded matrix defined by $\Phi^1$, $\Gamma^1$, $\Phi^2$, and $\Gamma^2$ if and only~if
  \[
    p^*_1(X) \leq b^*_2(X) \ \ \ \hbox{and}\ \ \ p^*_2(X) \leq b^*_1(X)
  \]
  hold for every subset $X \subseteq S_{m,n}$, where $p^*_1$, $b^*_1$, $p^*_2$, and $b^*_2$ are the functions defined in~\eqref{eq:bi*pi*}.
  For such matrices, the maximum possible total entry sum is $\min \{b^*_1(S_{m,n}), b^*_2(S_{m,n})\}$, and the minimum total entry sum is $\max \{p^*_1(S_{m,n}), p^*_2(S_{m,n})\}$.
  \FBOX
\end{corollary}

In the special case where $\alpha = -\infty$ and $\beta = +\infty$, Theorem~\ref{thm:general1} states that there exists an $(f, g)$-entry-bounded PBM if and only if~\eqref{eq:gen1a} and~\eqref{eq:gen1b} hold for every $X_1, X_2 \subseteq S_{m,n}$.

Our next goal is to show that this condition can be further simplified under the additional assumption that the last columns of $\Phi^1$ and $\Gamma^1$ are equal, and likewise, the last rows of $\Phi^2$ and $\Gamma^2$ are equal.
We refer to such a matrix as a {\bf strict prefix-bounded matrix}.
ASMs, by definition, fall into this stricter class.
In this case, the polyhedron $Q^1$ is the direct sum of $m$ base-polyhedra defined on the rows of $S_{m,n}$, while the polyhedron $Q^2$ is the direct sum of $n$ base-polyhedra defined on the columns of $S_{m,n}$, thus both $Q^1$ and $Q^2$ are base-polyhedra, as well.
By applying Corollary~\ref{cor:Bresz} to $Q_1$ and $Q_2$, we obtain the following characterization.
\begin{corollary}\label{cor:strictPBM}
  Assume that the last column of $\Phi^1$ is equal to that of $\Gamma^1$, and the last row of $\Phi^2$ is equal to that of $\Gamma^2$.
  Furthermore, assume that the total sum of the entries in the last column of $\Phi^1$ (and thus that of $\Gamma^1$) equals the total sum of the entries in the last row of $\Phi^2$ (and thus that of $\Gamma^2$); and let $H$ denote this common value.
  Then, for these prefix bounds $\Phi^1 \leq \Gamma^1$ and $\Phi^2 \leq \Gamma^2$, there exists a strict prefix-bounded matrix whose entries are bounded by $f$ and $g$ if and only if
  \begin{equation}\label{eq:pbmstrict}
    H\leq b^*_1(X') + b^*_2(X'') + \widetilde g(\ol{X'} \cap \ol{X''}) - \widetilde f (X' \cap X'')
  \end{equation}
  holds for every $X',X'' \subseteq S_{m,n}$, where $b^*_1$ and $b^*_2$ are the functions defined in~\eqref{eq:bi*pi*}.
  If there are no bounds imposed on the entries (that is, $f \equiv -\infty$ and $g \equiv \infty$), then the condition simplifies to
  \begin{equation}\label{eq:pbmstrict2}
    H\leq b^*_1(X) + b^*_2(\ol X)
  \end{equation}
  for every $X \subseteq S_{m,n}$.
\end{corollary}


\begin{remark}
  The concept of prefix-bounded matrices naturally extends to even more general settings.
  In addition to imposing bounds on horizontal- and vertical-prefix sums, one may also impose lower and upper bounds on the total sum of entries across the first $i$ rows (for each $i \in [m]$) and the first $j$ columns (for each $j \in [n]$).
  It is apparent that these additional requirements can be incorporated into the framework of laminar g-polymatroids.
  Consequently, the results established above extend naturally to this broader setting.
  In fact, the previously discussed constraint on the total sum of entries can be seen as a special case of this generalization.
  $\bullet$
\end{remark}

\begin{remark}
  The results of this section on prefix-bounded matrices naturally extend to bipartite graphs.
  Suppose we are given a bipartite graph $G = (V_1, V_2; E)$, where the edges incident to each vertex $v \in V_1 \cup V_2$ are endowed with a total order.
  Let $\Phi^1 \leq \Gamma^1$ and $\Phi^2 \leq \Gamma^2$ be functions defined on $E$, serving as lower and upper bounds.
  We are interested in the polyhedron of vectors $x \in \R^E$ satisfying the following constraints: for every edge $e = v_1v_2 \in E$, the sum of the $x$-values over the prefix of edges incident to $v_1$ ending with $e$ is between $\Phi^1(e)$ and $\Gamma^1(e)$; similarly, the sum of $x$-values over the prefix of edges incident to $v_2$ ending with $e$ lies between $\Phi^2(e)$ and $\Gamma^2(e)$.
  This polyhedron is the intersection of two prefix-bounded g-polymatroids, one for each vertex class, and as such, all previous results extend to this setting with minimal modification.
  $\bullet$
\end{remark}

\begin{remark}\label{rem:NPc}
  It is natural to consider additional constraints on the suffixes of the matrix, in addition to entry and prefix bounds.
  However, this generalization makes the problem hard, as it contains the NP-complete simultaneous matching problem~\cite{kutz2008simultaneous} as a special case.
  In fact, deciding the existence of an element-, prefix-, and suffix-bounded $m \times n$ integer matrix is NP-complete, even when suffix-bounds are imposed only on the rows.
  Another natural variation is to impose bounds on horizontal and vertical segments of length $d$ instead of prefixes, where the integer $d$ is a part of the input.
  Deciding whether an $m \times n$ integer matrix exists that satisfies entry bounds and, within each row and column, uniform bounds on the segments of length $d$ is also NP-complete, since it contains the perfect distance matching problem~\cite{madarasi2021matchings} as a special case.
  The problem remains hard even when $d$ equals the number of rows, all lower bounds are uniformly $0$, and the upper bounds are either $0$ or $1$ depending on the row and columns.
  $\bullet$
\end{remark}

\begin{remark}\label{rem:uniform_segment}
  We now propose a generalization of ASMs distinct from PBMs, which we refer to as {\bf segment-bounded matrices}.
  Within each row and column, suppose we are given uniform lower and upper bounds (depending on the row or column) on the sums over all segments, together with lower and upper bounds on the total sum of the row or column, as well as on the entries.
  When the upper bound for the segments is uniformly $+1$, the lower bound is $+1$ for the total sum of each row and column, and $-1$ for every other segment, we obtain ASMs as a special case.
  By Remark~\ref{rem:single_row_suffment} and Proposition~\ref{prop:box_plank_intersection}, the convex hull of (integer) segment-bounded matrices is the intersection of two g-polymatroids, and hence the polymatroidal approach works in this case.
  To the best of our knowledge, deciding the existence of a segment-bounded matrix necessarily involves g-polymatroids; unlike PBMs, this problem does not seem to be tractable via network-flow techniques.
  In contrast, the problem becomes NP-complete if the segment bounds are not uniform per row and column, see Remark~\ref{rem:NPc}.
  $\bullet$
\end{remark}

\section{Special classes of prefix-bounded matrices}\label{sec:special-PBM}

In this section, we explore how the previously established results on prefix-bounded matrices apply to certain well-known matrix classes.

\subsection{Alternating sign matrices}

Observe that ASMs are exactly the $n \times n$ strict PBMs in the special case when the lower bound $\Phi^1$ is identically $+1$ in the last column and $0$ elsewhere; the lower bound $\Phi^2$ is identically $+1$ in the last row and $0$ elsewhere; and the upper bounds $\Gamma^1$ and $\Gamma^2$ are identically $+1$.

Recall the definitions of separated systems of horizontal and vertical segments from Section~\ref{sec:pb-g-poli}.
For a given subset $X \subseteq S_{m,n}$, let $\mathcal{I}^{\mathrm{ho}}(X)$ and $\mathcal{I}^{\mathrm{ve}}(X)$ denote the system of the maximal horizontal and vertical segments in $X$, respectively.
We introduce the notation
\[
  \begin{cases}
    & \sigma_1(X) \coloneqq \vert {\cal I}^{\rm ho}(X)\vert,\\
    & \sigma_2(X) \coloneqq \vert {\cal I}^{\rm ve}(X)\vert.
  \end{cases}
\]
Clearly, for any subset $X \subseteq S_{m,n}$ formed as the union of a system $\mathcal{I}$ of separated horizontal or vertical segments, we have $\mathcal{I} = \mathcal{I}^{\mathrm{ho}}(X)$ or $\mathcal{I} = \mathcal{I}^{\mathrm{ve}}(X)$, respectively.

\begin{theorem}\label{thm:AS-feas}
  Let $f\leq g$ be lower and upper bounds on $S_{n,n}$.
  The following three conditions are pairwise equivalent.
  \medskip

  \noindent {\bf (A)} There exists an alternating sign matrix whose entries lie between $f$ and $g$.

  \medskip

  \noindent {\bf (B1)}
  \begin{equation}\label{eq:nkgf1b}
    \vert {\cal I}^1\vert + \vert {\cal I}^2\vert \geq n + \widetilde f (I_2)- \widetilde g(I_0)
  \end{equation}
  holds for every system ${\cal I}^1$ of separated horizontal segments and for every system ${\cal I}^2$ of separated vertical segments, where $I_2$ denotes the set of those positions of $S_{n,n}$ which are covered by both ${\cal I}^1$ and ${\cal I}^2$, while $I_0$ denotes the set of those positions which are covered by neither of them.
  \medskip

  \noindent {\bf (B2)}
  \begin{equation}\label{eq:nkgf1}
    \sigma_1(X') + \sigma_2(X'') \geq n + \widetilde f (X' \cap X'')- \widetilde g(\ol{X'} \cap \ol{X''})
  \end{equation}
  holds for every $X',X'' \subseteq S_{n,n}$.
\end{theorem}
\begin{proof}
  The equivalence of~(B1) and~(B2) is immediate by the definitions of $\sigma_1$, $\sigma_2$, and by the relationships between systems of segment and subsets of $S_{n,n}$.
  To prove the necessity of~\eqref{eq:nkgf1}, assume that there exists an ASM $z$ satisfying the bounds $f \leq z \leq g$, and let $z(s)$ denote the entry of this matrix at position $s$.
  Then, by the definition of ASMs, $\widetilde z(I)\leq 1$ holds for each segment $I$.
  Therefore, for every $X', X'' \subseteq S_{n,n}$, we obtain
  \[
    \begin{cases}
      & \widetilde z(X') = \sum\ [\widetilde z(I) : I \in {\cal I}^{\rm ho}(X')] \leq \sigma_1(X'),\\
      & \widetilde z(X'') = \sum\ [\widetilde z(I) : I \in {\cal I}^{\rm ve}(X'')] \leq \sigma_2(X''),
    \end{cases}
  \]
  from which
  \begin{align*}
    \begin{split}
      n &= \widetilde z(S_{n,n}) = \widetilde z(X') + \widetilde z(X'') - \widetilde z(X' \cap X'') + \widetilde z(\ol{X'} \cap \ol{X''})\\
        &\leq \sigma_1(X') + \sigma_2(X'') - \widetilde f(X' \cap X'' ) + \widetilde g(\ol{X'} \cap \ol{X''})
    \end{split}
  \end{align*}
  follows, which proves the necessity of~\eqref{eq:nkgf1}.

  The sufficiency of the condition will be derived from Corollary~\ref{cor:strictPBM}.
  Consider the case when $m = n$, the lower bound $\Phi^1$ is identically $+1$ in the last column and $0$ elsewhere; $\Phi^2$ is identically $+1$ in the last row and $0$ elsewhere; and the upper bounds $\Gamma^1$ and $\Gamma^2$ are identically $+1$.
  Recall that the class of strict PBMs defined by these functions coincides precisely with the class of ASMs.

  Now the values of $b^*_1$ on the horizontal segments and the values of $b^*_2$ on the vertical segments are uniformly $+1$, and $H = n$.
  Thus, the condition in~\eqref{eq:pbmstrict} from Corollary~\ref{cor:strictPBM} reduces exactly to~\eqref{eq:nkgf1}, completing the proof.
\end{proof}

For ease of application, we reformulate Theorem~\ref{thm:AS-feas} in an equivalent but slightly more tangible form.
Let ${\cal S} \coloneqq \{S_0, S_{+1}, S_{-1}, S_+, S_-, S_F\}$ be a partition of the ground set $S_{n,n}$, where each part may be empty.
We say that a matrix is {\bf $\boldsymbol{{\cal S}}$-compatible} if each position in $S_0$, $S_{+1}$, and $S_{-1}$ is assigned the fixed value $0$, $+1$, and $-1$, respectively; each position in $S_+$ is assigned a non-negative value, and each position in $S_-$ is assigned a non-positive value.
No restrictions are imposed on the values of positions in $S_F$.
We refer to the positions in the set $S_P \coloneqq S_0 \cup S_{+1} \cup S_{-1}$ as the {\bf prescribed} positions, and those in $S_F$ as the {\bf free} positions.
Our goal is to characterize those partitions $\cal S$ for which there exists an $\cal S$-compatible ASM.

Note that a system ${\cal I}^1$ of separated horizontal segments is not necessarily disjoint from a system ${\cal I}^2$ of separated vertical segments, since it may happen that a one-element segment $\{s\}$ appears in both systems.
To accommodate this, we adopt the convention that the union of the two systems is denoted by ${\cal I} \coloneqq {\cal I}^1 \uplus {\cal I}^2$, where one-element segments occurring in both systems are included with multiplicity two.
That is, ${\cal I}$ is formally treated as a family of segments, rather than a system of segments.
Accordingly, when we refer to a separated family of (horizontal and vertical) segments, we allow a one-element segment to occur in two copies.

We say that the family ${\cal I} \coloneqq {\cal I}^1 \uplus {\cal I}^2$ of segments is {\bf $\boldsymbol{{\cal S}}$-feasible} if
\begin{equation}\label{eq:megenga}
  I_0 \ \subseteq \ S_0 \cup S_{-1} \cup S_{-} \ \ \ \hbox{and}\ \ \ I_2 \ \subseteq \ S_0 \cup S_{+1} \cup S_{+},
\end{equation}
where $I_0$ denotes the set of positions in $S_{n,n}$ which are not covered by any of these two systems, while $I_2$ is the set of positions which are covered by both ${\cal I}^1$ and ${\cal I}^2$.
Note that ${\cal S}$-feasibility is equivalent to requiring that $\cal I$ covers the positions in $S_F$ exactly once, the positions in $S_- \cup S_{-1}$ at most once, and the positions in $S_+ \cup S_{+1}$ at least once.

\begin{theorem}\label{thm:AS-feas2}
  There exists an $\cal S$-compatible alternating sign matrix if and only if
  \begin{equation}\label{eq:nkgf2}
    \vert {\cal I}\vert \geq n + \vert S_{-1} \cap I_0\vert + \vert
    S_{+1}\cap I_2\vert
  \end{equation}
  holds for every $\cal S$-feasible family ${\cal I} \coloneqq {\cal I}^1 \uplus {\cal I}^2$ of segments, where $I_0$ and $I_2$ denote the set of positions not covered by $\cal I$ and covered twice by $\cal I$, respectively.
\end{theorem}
\begin{proof}
  The $\cal S$-compatibility of an ASM can be encoded with the help of lower and upper bound functions $f$ and $g$ defined on the ground set $S_{n,n}$ as follows:
  \[
    f(s) \coloneqq
    \begin{cases}
      0 & \ \ \hbox{if}\ \ \ s\in S_0 \cup S_+,\\
      +1 & \ \ \hbox{if}\ \ \ s\in S_{+1},\\
      -\infty & \ \ \hbox{if}\ \ \ s \in S_{-1} \cup S_-\cup S_F,\\
    \end{cases}
    \qquad
    g(s) \coloneqq
    \begin{cases}
      0 & \ \ \hbox{if}\ \ \ s\in S_0 \cup S_-,\\
      -1 & \ \ \hbox{if}\ \ \ s\in S_{-1},\\
      +\infty & \ \ \hbox{if}\ \ \ s \in S_{+1} \cup S_+\cup S_F.\\
    \end{cases}
  \]

  Note that $f(s)$ is finite if and only if $s\in S_0 \cup S_{+1} \cup S_+$, so the second condition in~\eqref{eq:megenga} is equivalent to requiring that $\widetilde f(I_2)$ is finite; in this case, $\widetilde f(I_2) = \vert S_{+1} \cap I_2 \vert$.
  Similarly, $g(s)$ is finite if and only if $s\in S_0 \cup S_{-1} \cup S_-$, and hence the first condition in~\eqref{eq:megenga} is equivalent to requiring that $\widetilde g(I_0)$ is finite; in this case, $\widetilde g(I_0) = - \vert S_{-1} \cap I_0 \vert$.

  Now we apply Theorem~\ref{thm:AS-feas}.
  Clearly,~\eqref{eq:nkgf1b} holds automatically when $\widetilde f(I_2) = -\infty$ or $\widetilde g(I_0) = +\infty$, thus, we only need to require the inequality in cases where both quantities are finite.
  In such cases, substituting the expressions for $\widetilde f(I_2)$ and $\widetilde g(I_0)$ into~\eqref{eq:nkgf1b} yields the condition in~\eqref{eq:nkgf2}.
\end{proof}

Let us now consider some notable special cases of Theorem~\ref{thm:AS-feas2}.
If $S_{-} \cup S_{+} = \emptyset$, meaning that only free and prescribed positions are present and no sign constraints are imposed, then the equivalent definition of ${\cal S}$-feasibility yields the following simplified form of Theorem~\ref{thm:AS-feas2}.
\begin{corollary}
  Let $\{S_0, S_{-1}, S_{+1}, S_F\}$ be a partition of $S_{n,n}$.
  Then there exists an alternating sign matrix in which each entry with position in $S_i$ is equal to $i$ for $i \in \{0, -1, +1\}$ if and only if~\eqref{eq:nkgf2} holds for every separated family ${\cal I}$ of segments that covers each position in $S_F$ exactly once, each position in $S_{-1}$ at most once, and each position in $S_{+1}$ at least once.
  \FBOX
\end{corollary}

If $S_{-1} \cup S_{+1} = \emptyset$, meaning that only sign constraints may be imposed and neither $+1$ nor $-1$ values are prescribed, then Theorem~\ref{thm:AS-feas2} reduces to the following simpler form.
\begin{corollary}\label{cor:subord}
  Let $\{S_0,S_-, S_+, S_F \}$ be a partition of the set $S_{n,n}$.
  There exists an alternating sign matrix which is $0$ on the positions in $S_0$, non-negative on the positions in $S_+$, and non-positive on the positions in $S_-$ if and only if $\vert {\cal I}\vert \geq n$ holds for every separated family ${\cal I}$ of segments which covers the positions in $S_F$ exactly once, the positions in $S_+$ at least once, end the positions in $S_-$ at most once.
  \FBOX
\end{corollary}

Brualdi and Dahl~\cite{brualdi2024frobenius} introduced the following relationship between matrices.
An $n \times n$ matrix $A$ is said to be a {\bf subordinate} of an $n \times n$ matrix $X$ if $A$ can be obtained from $X$ by replacing some non-zero entries with zeros in $X$.
The problem considered in~\cite{brualdi2024frobenius} is to characterize those matrices $X$ which have an ASM subordinate.
Theorem~4.1 and~4.3 in their paper provide necessary and sufficient conditions for the existence of such a subordinate in two special cases.
Our next result extends their findings by providing a characterization that applies to the general case.
\begin{corollary}\label{cor:alarendelt}
  A $(0,\pm 1)$-valued $n \times n$ matrix $X$ has an alternating sign matrix subordinate if and only if it is not possible to cover the $+1$ entries of $X$ using fewer than $n$ separated segments in such a way that each $-1$ entry is covered by at most one segment.
\end{corollary}
\begin{proof}
  Finding a subordinate ASM of $X$ is the same as finding an ${\cal S}$-compatible ASM for the partition where $S_0$ is the set of positions where $X$ is $0$, $S_+$ is the set of positions where $X$ is $+1$, and $S_-$ is the set of positions where $X$ is $-1$, while $S_F, S_{-1}$, and $S_{+1}$ are empty.
  The corollary then follows directly by applying Corollary~\ref{cor:subord} to this setting.
\end{proof}

\begin{remark}
  In their recent paper~\cite{brualdi2024frobenius}, Brualdi and Dahl posed the problem (Question 4.5) of maximizing the number of $-1$ entries over the subordinate ASMs of a given matrix $X$.
  This is clearly equivalent to maximizing the number of $+1$ entries.
  As we have seen earlier, the convex hull of $\cal S$-feasible ASMs is the intersection of two laminar base-polyhedra, and therefore, there exists a strongly polynomial algorithm for computing a maximum $w$-weight $\cal S$-feasible ASM for any weight function $w : S_{n,n}\rightarrow \R$.
  In particular, if we define $w$ to be $+1$ on the positions where $X$ is $+1$ and $0$ elsewhere, then a maximum $w$-weight subordinate ASM of $X$ has the maximum number of $+1$ entries.
  Moreover, since the underlying base-polyhedra are laminar, we show in Section~\ref{sec:circulation} that the relevant algorithmic problems can be addressed efficiently using standard network-flow algorithms.
  $\bullet$
\end{remark}

\begin{remark}
  Brualdi and Kim~\cite{brualdi2015completions} established, by a self-contained inductive proof, that any so-called ``bordered-permutation'' $(0,-1)$-valued $n \times n$ matrix can be extended to an ASM by turning certain $0$ entries into $+1$.
  Our framework is well-suited to handle such extension problems; in particular, let $S_{-1}$ denote the set of positions where $A$ is $-1$, and let $S_+$ denote the set of positions where $A$ is $0$.
  Then Theorem~\ref{thm:AS-feas2} provides a necessary and sufficient condition for the existence of such an extension for any $(0,-1)$-valued starting matrix.
  One can verify that this condition is satisfied in the special case of bordered-permutation matrices, and thus, the result of Brualdi and Kim can be derived as a corollary.
  However, this derivation is not more straightforward than the original inductive proof due to Brualdi and Kim.

  The melody of this relationship between a general result and its special case is similar to the situation when we want to prove the following ``theorem''.
  {\it Every 2-regular bipartite graph has a perfect matching}.
  Although this fact can be deduced from Hall's theorem, verifying the Hall-condition for such a graph is not at all simpler than proving the existence of a perfect matching directly.
  $\bullet$
\end{remark}

\subsection{Partial alternating sign matrices}

Recall that a partial alternating sign matrix (PASM), first introduced and studied by Heuer and Striker~\cite{heuer2022partial}, is an $m \times n$ prefix-bounded matrix defined by the identically~$0$ bounding functions $\Phi^1$ and $\Phi^2$ and the identically $+1$ functions $\Gamma^1$ and $\Gamma^2$.
For this special class of PBMs, the strong bordering pair can be derived directly from~\eqref{eq:bi*pi*}, allowing us to apply Theorem~\ref{thm:general1} and obtain the following result.

Recall that $\sigma_1(X)$ and $\sigma_2(X)$ denote the number of maximal horizontal and vertical segments in $X$, respectively.
Let $\sigma^{\operatorname{se}}_1(X)$ and $\sigma^{\operatorname{se}}_2(X)$ denote the number of maximal proper horizontal and vertical segments (that is, segments that are neither prefixes nor suffixes) in $X$, respectively, where the superscript ``$\operatorname{se}$'' stands for segment.
Similarly, let $\sigma^{\operatorname{su}}_1(X)$ and $\sigma^{\operatorname{su}}_2(X)$ denote the number of maximal proper horizontal and vertical suffixes (that is, suffixes that are not full rows and columns) in $X$, respectively, where the superscript ``$\operatorname{su}$'' stands for suffix.

\begin{theorem}\label{thm:PASM}
  There exists a partial alternating sign matrix of size $m \times n$ whose entries lie between $f$ and $g$, and whose total sum lies between $\alpha$ and $\beta$ if and only if the following inequalities hold for every $X_1, X_2 \subseteq S_{m,n}$:
  \begin{gather}
    p^*_1(X_1) + \widetilde f(X_2-X_1) \leq b^*_2(X_2) + \widetilde g(X_1-X_2),\label{eq:pasma}\\
    p^*_2(X_2) + \widetilde f(X_1-X_2) \leq b^*_1(X_1) + \widetilde g(X_2-X_1),\label{eq:pasmb}\\
    b^*_1(X_1) + b^*_2(X_2) + \widetilde g(\ol{X_1} \cap \ol{X_2}) - \widetilde f(X_1\cap X_2) \geq \alpha,\label{eq:pasmalfa}\\
    p^*_1(X_1) + p^*_2(X_2) + \widetilde f(\ol{X_1} \cap \ol{X_2}) - \widetilde g(X_1\cap X_2) \leq \beta,\label{eq:pasmbeta}
  \end{gather}
  where, for $i = 1,2$, the functions $p^*_i$ and $b^*_i$ are defined for $X \subseteq S_{m,n}$ as
  \begin{equation*}
    \begin{aligned}
      &\begin{cases}
        & \hbox{$p^*_i(X) \coloneqq -\sigma^{\operatorname{se}}_i(X) - \sigma^{\operatorname{su}}_i(X)$,}\\
        & \hbox{$b^*_i(X) \coloneqq \sigma_i(X).$}\\
      \end{cases}
    \end{aligned}
  \end{equation*}
  \FBOX
\end{theorem}

In the special case where $\alpha = -\infty$ and $\beta = +\infty$, Theorem~\ref{thm:PASM} states that there exists an $(f, g)$-entry-bounded PASM if and only if~\eqref{eq:pasma} and~\eqref{eq:pasmb} hold for every $X_1, X_2 \subseteq S_{m,n}$.
Finally, we note that both Corollaries~\ref{cor:maxpkm} and~\ref{cor:nofg} remain valid in the special case of PASMs.

\medskip
Both Behrend and Knight~\cite{behrend2007higher} and Striker~\cite{striker2007alternating, striker2009alternating} proved that the integral elements of the polytope $P$ of ASMs are in fact vertices of $P$.
In a forthcoming work, we will show the following general structural result, which immediately implies the former property of the polytope of ASMs.
Moreover, it also implies the analogous property of the polytope of PASMs. 
\begin{theorem}
  Let $P$ denote the convex hull of all $(0,\pm1)$-valued $m \times n$ matrices in which, for each row, the non-zero entries alternate in sign and the sign of either the first or the last non-zero entry is specified.
  Then, for any polyhedron $Q$, every integral element of $P \cap Q$ is a vertex of $P \cap Q$.
  \FBOX
\end{theorem}

\subsection{Alternating sign matrices with border constraints}\label{sec:borderedASM}

Brualdi and Kim~\cite{brualdi2015generalization, brualdi2021generalized, kim2024bruhat} introduced and investigated the following natural generalization of ASMs.
A $(0,\pm1)$-valued $m \times n$ matrix $A$ is a {\bf $\boldsymbol{(u,u'|v,v')}$-ASM} for border vectors $u,u'\in \{-1,+1\}^n$ and $v, v'\in \{-1,+1\}^m$ if the non-zero entries of row $i$ preceded by $v_i$ and succeeded by $v'_i$ and the non-zero entries of column $j$ preceded by $u_j$ and succeeded by $u'_j$ alternate in sign.
In other words, the non-zero entries of each row and column of the $(m+2)\times(n+2)$ matrix in Figure~\ref{fig:bordered} alternate in sign, except for the first and last rows and columns.
\begin{figure}[h]
  \centering
  \renewcommand{\arraystretch}{1.2}
  \renewcommand{\tabcolsep}{0.2cm}
  \begin{NiceTabular}{|wc{0.18cm}||c||wc{0.18cm}|}
    \hline
    $0$ & $u$ & $0$ \\
    \hline
    \hline
    $v$ & \Gape[0.8cm][0.8cm]{\hspace{1cm}$A$\hspace{1cm}} & $v'$ \\
    \hline
    \hline
    $0$ & $u'$ & $0$ \\
    \hline
  \end{NiceTabular}
  \caption{Illustration of the $(m+2)\times(n+2)$ matrix in the definition of $(u,u'|v,v')$-ASMs.}\label{fig:bordered}
\end{figure}

In this section, we derive a necessary and sufficient condition for the existence of a $(u,u'|v,v')$-ASM, which is another form of the condition presented in~\cite{brualdi2015generalization}.
In fact, our approach extends to a more general setting as we will explain below.

We embed this problem into the framework of prefix-bounded matrices.
We define the {\bf border type} of each row or column as $(-,-)$, $(+,+)$, $(-,+)$, or $(+,-)$, based on the sign of the corresponding entries of the border vectors $v,v'$ or $u,u'$, respectively.
For example, the border type of row $i$ is $(+,-)$ if $v_i = +1$ and $v'_i = -1$.

For a positive integer $k$ and an integer $\ell \in [k]$, define the values
\begin{equation*}
  \begin{aligned}
    &\begin{cases}
      \gamma^{--}(\ell) \coloneqq +1  & \text{for } \ell \in [k-1], \\
      \phi^{--}(\ell)   \coloneqq  0  & \text{for } \ell \in [k-1], \\
      \mathrlap{\gamma^{--}(k) \coloneqq \phi^{--}(k) \coloneqq +1,}\hphantom{\gamma^{++}(k) \coloneqq \phi^{++}(k) \coloneqq +1\quad}
    \end{cases} \\
    &\begin{cases}
      \gamma^{++}(\ell) \coloneqq  0 & \text{for } \ell \in [k-1], \\
      \phi^{++}(\ell)   \coloneqq -1 & \text{for } \ell \in [k-1], \\
      \mathrlap{\gamma^{++}(k) \coloneqq \phi^{++}(k) \coloneqq -1,}\hphantom{\gamma^{++}(k) \coloneqq \phi^{++}(k) \coloneqq +1\quad}
    \end{cases} \\
    &\begin{cases}
      \gamma^{-+}(\ell) \coloneqq +1 & \text{for } \ell \in [k-1], \\
      \phi^{-+}(\ell)   \coloneqq  0 & \text{for } \ell \in [k-1], \\
      \mathrlap{\gamma^{-+}(k) \coloneqq \phi^{-+}(k) \coloneqq 0,}\hphantom{\gamma^{++}(k) \coloneqq \phi^{++}(k) \coloneqq +1\quad}
    \end{cases} \\
    &\begin{cases}
      \gamma^{+-}(\ell) \coloneqq  0 & \text{for } \ell \in [k-1], \\
      \phi^{+-}(\ell)   \coloneqq -1 & \text{for } \ell \in [k-1], \\
      \mathrlap{\gamma^{+-}(k) \coloneqq \phi^{+-}(k) \coloneqq 0.}\hphantom{\gamma^{++}(k) \coloneqq \phi^{++}(k) \coloneqq +1\quad}
    \end{cases}
  \end{aligned}
\end{equation*}
\medskip

Let $k \coloneqq n$, and define the bounding functions for the horizontal prefixes as $\Phi^1(i,j) \coloneqq \phi^{s_1s_2}(j)$ and $\Gamma^1(i,j) \coloneqq \gamma^{s_1s_2}(j)$, where $(s_1,s_2)$ is the border type of row $i$.
Similarly, let  $k \coloneqq m$, and define the bounding functions for the vertical prefixes as $\Phi^2(i,j) \coloneqq \phi^{s_1s_2}(i)$ and $\Gamma^2(i,j) \coloneqq \gamma^{s_1s_2}(i)$, where $(s_1,s_2)$ is the border type of column $j$.

Now we derive a characterization for the existence of $(u,u'|v,v')$-ASMs.
Recall that $\sigma^{\operatorname{se}}_1(Z)$ and $\sigma^{\operatorname{se}}_2(Z)$ denote the number of maximal proper horizontal and vertical segments in $Z$, respectively, while $\sigma^{\operatorname{su}}_1(Z)$ and $\sigma^{\operatorname{su}}_2(Z)$ denote the number of maximal proper horizontal and vertical suffixes in $Z$.
Similarly, we define $\sigma^{\operatorname{pr}}_1(Z)$ and $\sigma^{\operatorname{pr}}_2(Z)$ to denote the maximal proper horizontal and vertical prefixes in $Z$, where the superscript ``$\operatorname{pr}$'' stands for prefix.
Similarly, $\sigma^{\operatorname{fu}}_1(Z)$ and $\sigma^{\operatorname{fu}}_2(Z)$ count the total number of full rows and columns in $Z$, with the superscript ``$\operatorname{fu}$'' standing for full.
By applying Corollary~\ref{cor:strictPBM} (for $f \equiv -\infty$ and $g \equiv \infty$) to this setting, we obtain the following result.
\begin{corollary}\label{cor:uuvv-ASM-existence}
  There exists a $(u,u'|v,v')$-ASM if and only if $b^*_1(S_{m,n}) = b^*_2(S_{m,n})$ and, for every $X \subseteq S_{m,n}$, the inequality
  \[
    H \leq b_1^*(X) + b_2^*(S_{m,n}-X)
  \]
  holds, where $H \coloneqq b^*_1(S_{m,n}) = b^*_2(S_{m,n})$, and for $i = 1, 2$,
  \[
    b_i^*(Z) \coloneqq \sum\ [b_i'(Z') : Z' \in {\cal C}^i(Z)].
  \]
  Here ${\cal C}^1(Z) \coloneqq \{ Z \cap P^1_{i,n} : i \in [m] \}$ and ${\cal C}^2(Z) \coloneqq \{ Z \cap P^2_{m,j} : j \in [n] \}$ denote the system of rows and columns of $Z$, respectively; and for $Z' \in {\cal C}^i(Z)$,
  \[
    b_i'(Z') \coloneqq
    \begin{cases}
      \sigma^{\operatorname{se}}_i(Z') + \sigma^{\operatorname{pr}}_i(Z') + \sigma^{\operatorname{su}}_i(Z') + \sigma^{\operatorname{fu}}_i(Z')    & \hbox{if $Z'$ is of border type $(-,-)$},\\
      \sigma^{\operatorname{se}}_i(Z') - \sigma^{\operatorname{fu}}_i(Z')    & \hbox{if $Z'$ is of border type $(+,+)$},\\
      \sigma^{\operatorname{se}}_i(Z') + \sigma^{\operatorname{pr}}_i(Z')    & \hbox{if $Z'$ is of border type $(-,+)$},\\
      \sigma^{\operatorname{se}}_i(Z') + \sigma^{\operatorname{su}}_i(Z')    & \hbox{if $Z'$ is of border type $(+,-)$},
    \end{cases}
  \]
  where the border type of a subset $Z'$ of a row (resp., column) is the border type of the row (resp., column) containing it.
  \FBOX
\end{corollary}

Corollary~\ref{cor:strictPBM} also provides necessary and sufficient condition for the existence of $(u,u'|v,v')$-ASMs subject to given lower and upper bounds on its entries.
Furthermore, the border vectors $u,u',v,v'$ may contain $0$ components, which can be handled as a projection of PBMs.

We emphasize that in Section~\ref{sec:circulation}, an alternative proof of Corollary~\ref{cor:strictPBM} will be presented using network-flow techniques.
This approach highlights that the existence problem for $(u,u'|v,v')$-ASMs is algorithmically tractable and can be solved efficiently using standard network-flow algorithms.

\begin{remark}
  Brualdi and Kim originally aimed to study the class of $(0, \pm 1)$-valued matrices whose non-zero entries in each row and column alternate in sign, with the first and last non-zero entries in each row and column prescribed.
  To encode these boundary prescriptions, they introduced the notion of $(u,u'|v,v')$-ASMs, where $u,u' \in \{\pm 1\}^m$ and $v,v' \in \{\pm 1\}^n$ specify the (negatives of the) desired first and last non-zero entries in each row and column, respectively.
  However, this framework does not faithfully represent the intended class of matrices, in particular, the resulting definition is stricter than the intended conditions.
  For instance, take $m = n = 1$ with $v_1 = v'_1 = +1$ and $u_1 = u'_1 = -1$.
  Clearly, no such $(u,u'|v,v')$-ASM exists, yet the $1 \times 1$ zero matrix satisfies the intended conditions.
  Consequently, Theorem 2.1 in the paper of Brualdi and Kim characterizes the existence of $(u,u'|v,v')$-ASMs but does not settle the original boundary-specified problem, which --- to the best of our knowledge --- remains open.
  Despite this mismatch, the $(u,u'|v,v')$-ASM framework yields a rich combinatorial structure with ties to the six-vertex model.
  Moreover, the results obtained in this setting admit natural generalizations beyond their original scope, as we shall see next.
  $\bullet$
\end{remark}

In what follows, we briefly outline the connection between ASMs --- and more generally $(u,u'|v,v')$-ASMs --- and the six-vertex model, also known as the square ice model~\cite{bressoud1999proofs}.
A {\bf state} of the six-vertex model is an $m \times n$ grid on $S_{m,n}$ where every vertex is incident to four directed arcs --- two vertical and two horizontal --- even along the boundary where extra boundary arcs are included, see Figure~\ref{fig:square_ice}.
The direction of the arcs is such that every vertex has exactly two incoming and two outgoing arcs.
This orientation constraint implies that the arcs incident to each vertex must be in one of the six states shown in Figures~\ref{fig:stateA}--\ref{fig:stateF}, called the {\bf state of the vertex}.

\begin{figure}[H]
  \centering
  \begin{subfigure}[t]{.16\textwidth}
    \centering
    \begin{tikzpicture}[scale = 0.9]
      \filldraw[] (0,0) circle (2pt);
      \draw[] (-1,0) -- (1,0);
      \draw[] (0,-1) -- (0,1);
      \draw[\niceArrow] (0,0) -- (0,0.6);
      \draw[\niceArrow] (0,0) -- (0,-0.6);
      \draw[\niceArrow] (-1,0) -- (-0.4,0);
      \draw[\niceArrow] (1,0) -- (0.4,0);
    \end{tikzpicture}
    \caption{}\label{fig:stateA}
  \end{subfigure}
  \hfill
  \begin{subfigure}[t]{.16\textwidth}
    \centering
    \begin{tikzpicture}[scale = 0.9]
      \filldraw[] (3,0) circle (2pt);
      \draw[] (3-1,0) -- (3+1,0);
      \draw[] (3,-1) -- (3,1);
      \draw[\niceArrow] (3,1) -- (3,0.4);
      \draw[\niceArrow] (3,-1) -- (3,-0.4);
      \draw[\niceArrow] (3,0) -- (2.4,0);
      \draw[\niceArrow] (3,0) -- (3.6,0);
    \end{tikzpicture}
    \caption{}\label{fig:stateB}
  \end{subfigure}
  \hfill
  \begin{subfigure}[t]{.16\textwidth}
    \centering
    \begin{tikzpicture}[scale = 0.9]
      \filldraw[] (6,0) circle (2pt);
      \draw[] (5,0) -- (7,0);
      \draw[] (6,-1) -- (6,1);
      \draw[\niceArrow] (6,0) -- (6,0.6);
      \draw[\niceArrow] (6,-1) -- (6,-0.4);
      \draw[\niceArrow] (5,0) -- (6-0.4,0);
      \draw[\niceArrow] (6,0) -- (6.6,0);
    \end{tikzpicture}
    \caption{}\label{fig:stateC}
  \end{subfigure}
  \hfill
  \begin{subfigure}[t]{.16\textwidth}
    \centering
    \begin{tikzpicture}[scale = 0.9]
      \filldraw[] (9,0) circle (2pt);
      \draw[] (8,0) -- (10,0);
      \draw[] (9,-1) -- (9,1);
      \draw[\niceArrow] (9,1) -- (9,0.4);
      \draw[\niceArrow] (9,0) -- (9,-0.6);
      \draw[\niceArrow] (10,0) -- (9.4,0);
      \draw[\niceArrow] (9,0) -- (8.4,0);
    \end{tikzpicture}
    \caption{}
  \end{subfigure}
  \hfill
  \begin{subfigure}[t]{.16\textwidth}
    \centering
    \begin{tikzpicture}[scale = 0.9]
      \filldraw[] (12,0) circle (2pt);
      \draw[] (11,0) -- (13,0);
      \draw[] (12,-1) -- (12,1);
      \draw[\niceArrow] (12,0) -- (12,0.6);
      \draw[\niceArrow] (12,-1) -- (12,-0.4);
      \draw[\niceArrow] (13,0) -- (12.4,0);
      \draw[\niceArrow] (12,0) -- (11.4,0);
    \end{tikzpicture}
    \caption{}
  \end{subfigure}
  \hfill
  \begin{subfigure}[t]{.16\textwidth}
    \centering
    \begin{tikzpicture}[scale = 0.9]
      \filldraw[] (15,0) circle (2pt);
      \draw[] (14,0) -- (16,0);
      \draw[] (15,-1) -- (15,1);
      \draw[\niceArrow] (15,1) -- (15,0.4);
      \draw[\niceArrow] (15,0) -- (15,-0.6);
      \draw[\niceArrow] (14,0) -- (14.6,0);
      \draw[\niceArrow] (15,0) -- (15.6,0);
    \end{tikzpicture}
    \caption{}\label{fig:stateF}
  \end{subfigure}
  \caption{The six possible states of the vertices.}
\end{figure}

There exists a well-known bijection between ASMs and the states of the $n\times n$ six-vertex model under {\bf domain wall boundary conditions}~\cite{elkies1992alternating2,kuperberg1996another}, which ensure that every extra horizontal arc on the border must point inward, and the extra vertical boundary arcs must point outward.
Given a state of the six-vertex model, the associated ASM $A = (a_{i,j})$ is defined by setting $a_{i,j} = +1$ if the vertex at position $(i,j)$ of the grid is in the state shown in Figure~\ref{fig:stateA}, $a_{i,j} = -1$ if it matches Figure~\ref{fig:stateB}, and $a_{i,j} = 0$ otherwise (i.e., if it corresponds to one of Figures~\ref{fig:stateC}-\ref{fig:stateF}).

This bijection extends naturally to $(u,u'|v,v')$-ASMs, which correspond to the states of the six-vertex model under {\bf general domain wall boundary conditions}, where the extra arcs on the boundary can be directed arbitrarily.
Specifically, an extra horizontal arc on the boundary pointing inward corresponds to a $-1$ component in the border vector $v$ (on the left side) or $v'$ (on the right side), while an outward-pointing arc corresponds to a $+1$ component.
Extra vertical boundary arcs are interpreted analogously, with inward and outward directions corresponding to $+1$ and $-1$ components, respectively, in $u$ or $u'$.
An example of a state of the six-vertex model with general domain wall boundary conditions can be seen in Figure~\ref{fig:square_ice} along with the corresponding $(u,u'|v,v')$-ASM.

\begin{figure}[H]
  \centering
  \begin{minipage}{0.45\textwidth}
    \begin{center}
      \begin{tikzpicture}
        \filldraw[black] (0,0) circle (2pt);
        \filldraw[black] (0,1) circle (2pt);
        \filldraw[black] (0,2) circle (2pt);
        \filldraw[black] (0,3) circle (2pt);
        \filldraw[black] (1,0) circle (2pt);
        \filldraw[black] (1,1) circle (2pt);
        \filldraw[black] (1,2) circle (2pt);
        \filldraw[black] (1,3) circle (2pt);
        \filldraw[black] (2,0) circle (2pt);
        \filldraw[black] (2,1) circle (2pt);
        \filldraw[black] (2,2) circle (2pt);
        \filldraw[black] (2,3) circle (2pt);
        \draw[] (-1,0) -- (3,0);
        \draw[] (-1,1) -- (3,1);
        \draw[] (-1,2) -- (3,2);
        \draw[] (-1,3) -- (3,3);
        \draw[] (0,-1) -- (0,4);
        \draw[] (1,-1) -- (1,4);
        \draw[] (2,-1) -- (2,4);
        \draw[\niceArrow] (0,0) -- (-0.6,0);
        \draw[\niceArrow] (-1,1) -- (-0.4,1);
        \draw[\niceArrow] (0,2) -- (-0.6,2);
        \draw[\niceArrow] (-1,3) -- (-0.4,3);
        \draw[\niceArrow] (0,-1) -- (0,-0.4);
        \draw[\niceArrow] (1,0) -- (1,-0.6);
        \draw[\niceArrow] (2,0) -- (2,-0.6);
        \draw[\niceArrow] (0,4) -- (0,3.4);
        \draw[\niceArrow] (1,3) -- (1,3.6);
        \draw[\niceArrow] (2,4) -- (2,3.4);
        \draw[\niceArrow] (0,0) -- (0.6,0);
        \draw[\niceArrow] (0,1) -- (0,0.4);
        \draw[\niceArrow] (1,0) -- (1.6,0);
        \draw[\niceArrow] (1,1) -- (1,0.4);
        \draw[\niceArrow] (3,0) -- (2.4,0);
        \draw[\niceArrow] (2,0) -- (2,0.6);
        \draw[\niceArrow] (0,1) -- (0.6,1);
        \draw[\niceArrow] (0,2) -- (0,1.4);
        \draw[\niceArrow] (2,1) -- (1.4,1);
        \draw[\niceArrow] (1,1) -- (1,1.6);
        \draw[\niceArrow] (2,1) -- (2.6,1);
        \draw[\niceArrow] (2,2) -- (2,1.4);
        \draw[\niceArrow] (1,2) -- (0.4,2);
        \draw[\niceArrow] (0,3) -- (0,2.4);
        \draw[\niceArrow] (1,2) -- (1.6,2);
        \draw[\niceArrow] (1,3) -- (1,2.4);
        \draw[\niceArrow] (2,2) -- (2.6,2);
        \draw[\niceArrow] (2,3) -- (2,2.4);
        \draw[\niceArrow] (3,3) -- (2.4,3);

        \draw[\niceArrow] (2,3) -- (1.4,3);
        \draw[\niceArrow] (0,3) -- (0.6,3);

      \end{tikzpicture}
    \end{center}
  \end{minipage}%
  \begin{minipage}{0.45\textwidth}
    \begin{center}
      \renewcommand{\arraystretch}{1.9}
      \renewcommand{\tabcolsep}{0.2cm}
      \begin{NiceTabular}{|c||ccc||c|}
        \hline
        $0$ & $+1$ & $-1$ & $+1$ & $0$\\
        \hline
        \hline
        $-1$ & $0$ & $+1$ & $0$ & $-1$\\
        $+1$ & $0$ & $-1$ & $0$ & $+1$\\
        $-1$ & $0$ & $+1$ & $-1$ & $+1$\\
        $+1$ & $-1$ & $0$ & $+1$ & $-1$\\
        \hline
        \hline
        $0$ & $+1$ & $-1$ & $-1$ & $0$\\
        \hline
      \end{NiceTabular}
    \end{center}
  \end{minipage}
  \caption{A state of $4 \times 3$ the six-vertex model and the corresponding $(u,u'|v,v')$-ASM with the bordering vectors.}\label{fig:square_ice}
\end{figure}

This bijection shows that the problem considered in Corollary~\ref{cor:uuvv-ASM-existence} is equivalent to deciding whether there exists a state of the six-vertex model in which the boundary arcs are oriented according to the border vectors $u$, $u'$, $v$, and $v'$.

\begin{remark}
  More generally, one can decide whether there exists a state of the six-vertex model in which the orientations of certain arcs are prescribed, using either PBMs or Hakimi's Orientation Lemma~\cite{hakimi1965degrees}.
  $\bullet$
\end{remark}

\section{Prefix-bounded matrices and feasible circulations}\label{sec:circulation}

Throughout this section, we consider PBMs defined by prefix bounds $\Phi^1$, $\Gamma^1$, $\Phi^2$, $\Gamma^2$, and require that the entries of the PBMs in question lie between $f$ and $g$, while the total sum of its entries lies between $\alpha$ and $\beta$.
The corresponding polyhedron is the intersection of two laminar g-polymatroids, and thus admits a description via a linear inequality system whose constraint matrix is the transpose of the incidence matrix of the union of two laminar families, see Remark~\ref{rem:2-lamin} and Corollary~\ref{cor:boundedPBMPolytope}.
Any such system can be obtained as the projection (restriction) of a feasible circulation polyhedron.

In what follows, we work out the details of this approach for the special case of PBMs.
In particular, we show how Hoffman's classical theorem on feasible circulations~\cite{hoffman1958some,micchelli2003selected} can be used to provide an alternative proof of Theorem~\ref{thm:general1}.
This proof also implies that standard network-flow algorithms can be used for finding the desired PBM with entry bounds and total-sum bounds, or finding the certificate stated in Theorem~\ref{thm:general1} for the non-existence of such a matrix.
Note that the minimum cost PBM problem is also tractable within this circulation framework.

\subsection{The circulation model}\label{sec:arammodell}
In this section, we define a feasible circulation network whose integer-valued solutions directly correspond to the $m \times n$ PBMs.
Let $D = (V, A)$ be a directed graph constructed as follows.
The vertex set $V$ is partitioned into two disjoint subsets, denoted by $V^1$ and $V^2$, defined as
\[
  V^1 \coloneqq \{v^1_{i,j} : (i,j) \in S_{m,n}\} + v^1_0 \ \ \ \hbox{and}\ \ \ V^2 \coloneqq \{v^2_{i,j} : (i,j) \in S_{m,n}\} + v^2_0.
\]

Now we define the arc set $A$ by constructing four classes of arcs.
First, within the vertex set $V^1$, for each $i\in [m]$ and $j\in [n-1]$, include an arc $a^1_{i,j} \coloneqq v^1_{i,j+1} v^1_{i,j}$.
Additionally, for each $i\in [m]$, add an arc $a^1_{i,n } \coloneqq v^1_0 v^1_{i,n}$.
Let
\[
  A^1 \coloneqq \{a^1_{i,j} : (i,j) \in S_{m,n}\}
\]
denote the set of all these arcs.
Next, within the vertex set $V^2$, for each $i\in [m-1]$ and $j\in [n]$, include an arc $a^2 _{i,j} \coloneqq v^2_{i,j} v^2_{i+1,j}$.
For each $j\in [n]$, add an arc $a^2_{m,j} \coloneqq v^2_{m,j}v^2_0$.
These arcs form the set
\[
  A^2 \coloneqq \{a^2_{i,j} : (i,j) \in S_{m,n}\}.
\]
To connect $V^1$ to $V^2$, for each $i \in [m]$ and $j \in [n]$, include an arc $a_{i,j} \coloneqq v^1_{i,j}v^2_{i,j}$; and let
\[
  N \coloneqq \{a_{i,j} : (i,j)\in S_{m,n} \}
\]
denote the set of these arcs.
Finally, insert a single arc connecting the two extra vertices
\[
  a_0 \coloneqq v^2_0v^1_0.
\]
Figure~\ref{fig:flow} shows this construction for the case $n = m = 4$.

We will associate the arcs in $A^1$ with the horizontal prefixes of the desired PBM, the arcs in $A^2$ with the vertical prefixes, the arcs in $N$ with the positions of the entries, and the arc $a_0$ with the total sum of the entries.

The arc set $A^1 \cup A^2 + a_0$ forms a (directed) spanning tree $T$ of $D$, whereas the arc set $N$ forms a matching between the vertices in $V^1-v^1_0$ and $V^2-v^2_0$.
We refer to the arcs in $T$ as {\bf tree arcs} and to those in $N$ as {\bf non-tree arcs}.
In Figure~\ref{fig:flow}, the green continuous and black dotted arrows represent the tree arcs and the non-tree arcs, respectively.

\begin{figure}[t]
  \centering
  \begin{tikzpicture}[xscale=.9]
    \tikzset{VertexStyle/.append style = {minimum size = 22pt,inner sep=0pt}}
    \SetVertexMath

    \def\hsep{4} 
    \def\vsep{1.25} 
    \def\hgain{1} 
    \def\height{5.75} 

    \Vertex[x=3*\hsep + 1*\hgain, y=\height - 1.5*\vsep, L=v^1_0]{s}
    \Vertex[x=-1 * \hgain + 1 * \hsep, y=-4*\vsep, L=v^2_0]{t}

    \foreach \i in {1,...,4} {
      \foreach \j in {1,...,3} {
        \pgfmathsetmacro{\x}{(\j - 1)*\hsep + (4 - \i) * \hgain}
        \pgfmathsetmacro{\y}{\height - (\i - 1)*\vsep}
        \Vertex[x=\x, y=\y, L=v^1_{\i,\j}]{s\i\j}
      }
    }

    \foreach \i in {1,...,4} {
      \foreach \j in {1,...,3} {
        \pgfmathsetmacro{\x}{(\j - 1)*\hsep + (4 - \i) * \hgain}
        \pgfmathsetmacro{\y}{-(\i - 1)*\vsep}
        \Vertex[x=\x, y=\y, L=v^2_{\i,\j}]{t\i\j}
      }
    }

    \foreach \j in {1,...,3} {
      \foreach \i [evaluate=\i as \k using int(\i+1)] in {1,...,3} {
        \draw[color=ForestGreen,\niceArrow] (t\i\j) -- (t\k\j);
      }
      \draw[color=ForestGreen,\niceArrow] (t4\j) -- (t);
    }

    \foreach \i in {1,...,4} {
      \foreach \j in {1,...,3} {
        \draw[draw=white, double distance=\pgflinewidth, ultra thick, shorten <=1pt, shorten >=1pt] (s\i\j) to (t\i\j);
        \draw[dotted, dash pattern=on 1pt off 2pt, \niceArrow] (s\i\j) -- (t\i\j);
      }
    }

    \foreach \i in {1,...,4} {
      \foreach \j [evaluate=\j as \k using int(\j+1)] in {1,...,2} {
        \draw[draw=white, double distance=\pgflinewidth, ultra thick, shorten <=1pt, shorten >=1pt] (s\i\k) to (s\i\j);
        \draw[color=ForestGreen,\niceArrow] (s\i\k) -- (s\i\j);
      }
    }

    \foreach \i in {1,...,4} {
      \draw[draw=white, double distance=\pgflinewidth, ultra thick, shorten <=1pt, shorten >=1pt] (s) to (s\i3);
    }
    \foreach \i in {1,...,4} {
      \draw[color=ForestGreen,\niceArrow] (s) -- (s\i3);
    }

    \draw[color=ForestGreen,\niceArrow] (t) to[out = 0, in = -90, looseness = 1.3] (s);
  \end{tikzpicture}
  \caption{The circulation network for a $4 \times 3$ prefix-bounded matrix.
    The vertical arcs correspond to the matrix entries.}\label{fig:flow}
\end{figure}

Observe that every one-way circuit in $D$ contains the arc $a_0$, that is, $D-a_0$ is acyclic.
Furthermore, for any non-tree arc $a$, the unique circuit $C_a$ in $T+a$, called the {\bf fundamental circuit} belonging to $a$, is a one-way circuit.
A basic property of circulations is that every circulation~$z : A \rightarrow \R$ on the arc set of $D$ can be written in the form ${z = \sum [ z(a) \chi(C_a) : a\in N]}$, where $\chi(C_a)$ denotes the characteristic vector of the fundamental circuit $C_a$, see Proposition~4.3.1 in~\cite{frank2011connections}.

The arcs $a_{i,1}, \dots, a_{i,j}$ are precisely the non-tree arcs whose fundamental circuits contain the tree arc $a^1_{i,j} \in A^1$, similarly, the arcs $a_{1,j},\dots, a_{i,j}$ are precisely the non-tree arcs whose fundamental circuits contain the tree arc $a^2_{i,j}\in A^2$, furthermore, the special tree arc $a_0$, as already mentioned, is contained in every fundamental circuit.
These observations immediately yield the following.
\begin{claim}\label{cl:zafan}
  For every circulation $z$ in $D$, the $z$-values of the tree arcs can be expressed in terms of the $z$-values of the non-tree arcs as follows:
  \[
    \begin{cases}
      & z(a^1_{i,j}) = z(a_{i,1}) + \cdots + z(a_{i,j}),\\
      & z(a^2_{i,j}) = z(a_{1,j}) + \cdots + z(a_{i,j}),\\
      & z(a_0) = \sum\ [ z(e) : e\in N].
    \end{cases}
  \]
  \FBOX
\end{claim}

Note that the formulas in Claim~\ref{cl:zafan} for the $z$-values of the tree arcs can be concisely expressed in matrix form as $z_T = M_Dz_N$, where $z_T$ and $z_N$ denote the restrictions of $z$ to the tree arcs and non-tree arcs, respectively, and $M_D$ is the network matrix associated with the (directed) spanning tree $T$ of the digraph $D$.

To encode the prefix, entry and total-sum bounds within our circulation network, we define lower and upper bound functions $l$ and $u$ on the arc set $A$.
These functions determine the admissible flow values on each arc in accordance with the bounds defining the desired PBM.
For each $(i,j) \in S_{m,n}$, the arc $a^1_{i,j} \in A^1$ corresponds to the horizontal prefix that consists of the first $j$ positions in row $i$ of the PBM.
Accordingly, the bounds for each arc $a^1_{i,j} \in A^1$ are defined as
\begin{equation}\label{eq:a1ij}
  l(a^1_{i,j}) \coloneqq \Phi^1 ({i,j}) \ \ \ \hbox{and}\ \ \ u(a^1_{i,j}) \coloneqq \Gamma^1({i,j}).
\end{equation}
Similarly, for each arc $a^2_{i,j} \in A^2$, which corresponds to the vertical prefix consisting of the first $i$ positions in column $j$, the bounds are defined by
\begin{equation}\label{eq:a2ij}
  l(a^2_{i,j}) \coloneqq \Phi^2 ({i,j}) \ \ \ \hbox{and}\ \ \ u(a^2_{i,j}) \coloneqq \Gamma^2({i,j}).
\end{equation}
For the arcs in the set $N$, each of which corresponds to an individual entry, the lower and upper bounds are given by $f$ and $g$.
That is,
\begin{equation}
  l(a_{i,j}) \coloneqq f(i,j) \ \ \ \hbox{and}\ \ \ u(a_{i,j}) \coloneqq g(i,j).
\end{equation}
Finally, for the special arc $a_0$, which represents the total sum of the entries, we define the bounds as 
\begin{equation}\label{eq:BoundsEnd}
  l(a_0) \coloneqq \alpha \ \ \ \hbox{and}\ \ \ u(a_0) \coloneqq \beta.
\end{equation}

The following theorem shows that there is a one-to-one correspondence between PBMs and integer-valued feasible circulations.
\begin{theorem}\label{thm:pbm-aram}
  Consider the circulation network defined above on the digraph $D$ with the lower and upper bounds $l$ and $u$.

  {\bf (A)} \ Let $z \coloneqq (z_N,z_T)$ be an integer-valued circulation in $D$ for which $l\leq z\leq u$ (in short, $z$ is $(l,u)$-bounded).
  Consider the $m \times n$ matrix $M_z$ for which $M_z(i,j) \coloneqq z_N(a_{i,j})$ for $(i,j) \in S_{m,n}$.
  Then $M_z$ is a matrix which is {\bf (a)} prefix-bounded with respect to horizontal-prefix bounds $\Phi^1$ and $\Gamma^1$ and vertical-prefix bounds $\Phi^2$ and $\Gamma^2$, {\bf (b)} entry-bounded by lower and upper bounds $f$ and $g$, and {\bf (c)} the total entry sum of $M_z$ is between $\alpha$ and $\beta$; that is, the matrix $M_z$ meets the requirements of Theorem~\ref{thm:general1}.

  {\bf (B)} \ Conversely, let $M$ be a prefix-bounded matrix with entry bounds and total-sum bounds.
  Define the vector $z_M$ on the arc set of $D$ as
  \begin{equation*}
    \begin{aligned}
      &z_M(a^1_{i,j}) \coloneqq M(i,1) + \dots + M(i,j) & \hbox{ for every $a_{i,j}^1\in A^1$,}\\
      &z_M(a^2_{i,j}) \coloneqq M(1,j) + \dots + M(i,j) & \hbox{ for every  $a_{i,j}^2\in A^2$,}\\
      &z_M(a_{i,j}) \coloneqq M(i,j) & \mathrlap{\hbox{ for every $ a_{i,j}\in N$,}}\hphantom{\hbox{ for every $a_{i,j}^1\in A^1$,}}\\   &z_M(a_0) \coloneqq \sum\ [ M(i,j) : (i,j) \in S_{m,n}].\\
    \end{aligned}
  \end{equation*}
  Then $z_M$ is an $(l,u)$-bounded integer-valued circulation in the digraph $D$.
\end{theorem}
\begin{proof}
  (A) (a) \ To prove the prefix-boundedness of $M_z$, consider the prefix $P^1_{i,j}$ consisting of the first $j$ positions of the $i$-th row.
  The sum of the entries at the positions in $P^1_{i,j}$ is
  \[
    M_z(i,1) + \dots + M_z(i,j) =  z(a_{i,1}) + \dots + z (a_{i,j}) = z(a^1_{i,j}),
  \]
  which follows by the definition of $M_z$ and from Claim~\ref{cl:zafan}.

  Analogously, the sum of the entries at the positions in $P^2_{i,j}$ is
  \[
    M_z(1,j) + \dots + M_z(i,j) = z(a_{1,j}) + \dots + z(a_{i,j}) =  z(a^2_{i,j}).
  \]

  By~\eqref{eq:a1ij}~and~\eqref{eq:a2ij}, the $(l,u)$-boundedness of $z$ implies that $M_z$ meets the lower and upper bounds $\Phi^1,\Gamma^1$ for horizontal and  $\Phi^2, \Gamma^2$ for vertical prefixes.

  \medskip

  (b) \ The constraint $l_N\leq z_N\leq u_N$ concerning non-tree arcs is equivalent to $f(i,j) \leq M_z(i,j) \leq g(i,j)$, that is, the entries of $M_z$ are indeed between $f$ and $g$.

  \medskip

  (c) \ The total sum of the entries of the matrix $M_z$ is given by $\widetilde z(N)$, which is equal to $z(a_0)$ by Claim~\ref{cl:zafan}.
  Since $z$ is $(l,u)$-bounded, we have $\alpha = l(a_0) \leq z(a_0) \leq u(a_0) = \beta$, and thus the total sum lies between the prescribed bounds.
  \medskip

  (B) \ It follows directly from the definition of $z_M$ and by the construction of $D$ that $z_M$ is a circulation.
  We now verify that $z_M$ is $(l,u)$-bounded, that is, $l\leq z_M\leq u$.
  By the definition of $z_M$, the horizontal-prefix bounds, the vertical-prefix bounds, the entry bounds, and the total-sum bounds for $M$ are equivalent to the $(l,u)$-boundedness of $z_M$ for the arcs in $A_1$, $A_2$, $N$, and $\{a_0\}$, respectively.
  Therefore, $z_M$ is $(l,u)$-bounded, which completes the proof.
\end{proof}

In the next section, this circulation framework will be used to prove Theorem~\ref{thm:general1} from Hoffman's classical theorem on feasible circulations.

\subsection{A proof of Theorem~\ref{thm:general1} via circulations}


To prove necessity, we show that if there exists a pair of subsets $X_1,X_2$ violating one of the inequalities \eqref{eq:gen1a}-\eqref{eq:gen1beta}, then there exists no desired PBM.
It is easy to see that the left-hand sides of~\eqref{eq:gen1a} and~\eqref{eq:gen1b} serve as lower bounds for the sum of the entries in $X_1 \cup X_2$, while the right-hand sides serve as corresponding upper bounds.
Consequently, if either inequality fails to hold, then the upper bound is strictly smaller than the lower bound, implying that no desired PBM exists.
In~\eqref{eq:gen1alfa}, the left-hand side is an upper bound for the sum of the entries in the matrix, thus, if this bound is smaller than $\alpha$, then no desired PBM exists.
Similarly, the left-hand side of~\eqref{eq:gen1beta} is a lower bound for the sum of the entries, thus, if this bound exceeds $\beta$, then again no desired PBM exists.

Our final goal is to show that the circulation framework above makes it possible to prove sufficiency from Hoffman's classical theorem on feasible circulations.
For a digraph $D = (V,A)$ and for a subset $W\subseteq V$, let $\Delta ^{\operatorname{in}}(W) \coloneqq \{uv \in A : u \notin W, v \in W \}$ denote the set of arcs entering $W$, and let $\Delta ^{\operatorname{out}}(W) \coloneqq \{uv \in A : u \in W, v \notin W \}$ denote the set of arcs leaving $W$.
For a function $x : A \rightarrow \R$ and for a subset $W\subseteq V$, we use the notations $\varrho _x(W) \coloneqq \widetilde x (\Delta ^{\operatorname{in}}(W))$ and $\delta _x(W) \coloneqq \widetilde x (\Delta ^{\operatorname{out}}(W))$ to denote the sum of $x$-values of the arcs entering and leaving $W$, respectively.
Given a function $x : A \rightarrow \R$ and a subset $F \subseteq A$, we write $x \vert F$ to denote the function which is equal to $x$ on $F$ and $0$ on $A-F$.
For example, for a subset $Z \subseteq V$, the expression $\varrho_{x \vert F}(Z)$ is the sum of $x$-values on the arcs in $F$ entering $Z$.

\begin{theorem}[Hoffman's circulation theorem~\cite{hoffman1958some, micchelli2003selected}]
  Let $l$ and $u$ be lower and upper bound functions on the arc set of a digraph $D = (V,A)$ for which $l\leq u$.
  There is a circulation $z$ for which $l\leq z\leq u$ (in short, $z$ is $(l,u)$-bounded) if and only if $\varrho _u(W) \geq \delta _l(W)$ holds for every subset $W\subseteq V$ of vertices.
  If, in addition, $l$ and $u$ are integer-valued, then the circulation $z$ may be chosen integer-valued.
  \FBOX
\end{theorem}

Now we are ready to derive Theorem~\ref{thm:general1} from Hoffman's theorem.
Consider the digraph $D = (V,A)$ introduced in Section~\ref{sec:arammodell} along with the lower and upper bounds $l$ and $u$ defined in~\eqref{eq:a1ij}-\eqref{eq:BoundsEnd}.
By Part (A) of Theorem~\ref{thm:pbm-aram}, any feasible integer-valued circulation --- when it exists --- determines an entry-bounded and prefix-bounded matrix whose total entry sum lies between $\alpha$ and $\beta$.

Now suppose that no $(l,u)$-bounded integer-valued circulation exists in $D$.
By Hoffman's theorem, there must exist a subset $W \subseteq V$ that violates the Hoffman-condition, that is,
\begin{equation}\label{eq:Hsertes}
  \varrho _u(W) - \delta _l(W) <0.
\end{equation}

Theorem~\ref{thm:general1} characterizes the existence of the desired PBM in terms of the functions $p^*_1$, $b^*_1$, $p^*_2$ and $b^*_2$, so we need to establish the connection between these functions and our circulation network.
Consider a segment $Z \subseteq S_{m,n}$ lying in row $i$, where $i \in [m]$.
Let $(i,h)$ and $({i,k})$ denote the first and last positions of $Z$, respectively.
By applying the definition of $(p^*_1,b^*_1)$ given in~\eqref{eq:bi*pi*} to the segment $Z$, 
we obtain from~\eqref{eq:a1ij} that
\begin{equation}\label{eq:b1*p1*}
  \begin{cases}
    & p^*_1(Z) = \Phi^1({i,k}) - \Gamma^1({i,h-1}) = l(a^1_{i,k}) - u(a^1_{i,h-1}),\\
    & b^*_1(Z) = \Gamma^1({i,k}) - \Phi^1({i,h-1}) = u(a^1_{i,k}) - l(a^1_{i,h-1}),
  \end{cases}
\end{equation}
with the conventions $l(a^1_{i,h-1}) \coloneqq 0$ and $u(a^1_{i,h-1}) \coloneqq 0$ when $h = 1$.
For a subset $X\subseteq S_{m,n}$, $p ^*_1(X)$ and $b ^*_1(X)$ are equal to the sum of the $p^*_1$ and $b^*_1$ values of the maximal horizontal segments in $X$, respectively.

Similarly, consider a segment $Z$ in the $j$-th column of $S_{m,n}$, where $j \in [n]$.
Let $({h,j})$ and $({k,j})$ denote the first and last positions of $Z$, respectively.
By~\eqref{eq:a2ij}, we have
\begin{equation}\label{eq:b2*p2*}
  \begin{cases}
    & p^*_2(Z) = \Phi^2({k,j}) - \Gamma^2({h-1,j}) = l(a^2_{k,j}) - u(a^2_{h-1,j}),\\
    & b^*_2(Z) = \Gamma^2({k,j}) - \Phi^2({h-1,j}) = u(a^2_{k,j}) - l(a^2_{h-1,j}),
  \end{cases}
\end{equation}
with the conventions $l(a^2_{h-1,j}) \coloneqq 0$ and $u(a^2_{h-1,j}) \coloneqq 0$ when $h = 1$.
For a subset $X\subseteq S_{m,n}$, $p ^*_2(X)$ and $b ^*_2(X)$ are equal to the sum of the $p^*_2$ and $b^*_2$ values of the maximal vertical segments in $X$, respectively.

To prove the theorem, we proceed by analyzing four distinct cases determined by the structure of the set $W$ that violates the Hoffman-condition.
These cases give rise to four types of inequalities, each corresponding to one of the necessary conditions stated in Theorem~\ref{thm:general1}.

\medskip

{\bf Case 1.} \ $ v^1_0,v^2_0 \in W$.
Let $W^1 \coloneqq V^1 -W$ and $W^2 \coloneqq V^2-W$.
Furthermore, let
\[
  X_1 \coloneqq \{ ( {i,j}) \in S_{m,n} : v^1_{i,j} \in W^1 \} \ \ \ \hbox{and}\ \ \
  X_2 \coloneqq \{ ({i,j}) \in S_{m,n} : v^2_{i,j} \in W^2 \}.
\]
Consider what~\eqref{eq:Hsertes} means in terms of $X_1$ and $X_2$.
It follows from~\eqref{eq:b1*p1*} and~\eqref{eq:b2*p2*} that the contributions of the arcs in $A^1$ and $A^2$ to $\varrho_u(W) - \delta_l(W)$ are
\begin{align*}
  \begin{split}
    &\varrho _{u\vert A^1}(W)- \delta _{l\vert A^1}(W) = \widetilde u(A^1 \cap \Delta^{\operatorname{out}}(W^1)) - \widetilde l(A^1 \cap \Delta^{\operatorname{in}}(W^1)) = - p^*_1(X_1)\\
    &\varrho _{u\vert A^2}(W) - \delta _{l\vert A^2}(W) = \widetilde u(A^2 \cap \Delta^{\operatorname{out}}(W^2)) - \widetilde l(A^2 \cap \Delta^{\operatorname{in}}(W^2)) = b^*_2(X_2).
  \end{split}
\end{align*}
The contribution of an arc $a_{i,j} = v^1_{i,j} v^2_{i,j} \in N$ to $\varrho_u(W)$ and $\delta_l(W)$ is $0$ if both endpoints of $a_{i,j}$ lie in $W$ or if neither does.
If $a_{i,j}$ enters $W$, that is, if $(i,j) \in X_1 - X_2$, then it contributes $g(i,j)$ to $\varrho_u(W)$ and $0$ to $\delta_l(W)$.
If $a_{i,j}$ leaves $W$, that is, if $(i,j) \in X_2 - X_1$, then it contributes $f(i,j)$ to $\delta_l(W)$ and $0$ to $\varrho_u(W)$.
Hence,
\[
  \varrho _{u\vert N}(W) - \delta _{l\vert N}(W) = \widetilde g(X_1-X_2) - \widetilde f(X_2-X_1).
\]
The contribution of the special arc $a_0$ is $0$, since both of its endpoints are in $W$.

Summing all contributions, we obtain that
\begin{align*}
  \begin{split}
    0 &> \varrho _u(W) - \delta _l(W) = [\varrho _{u\vert A^1} (W) + \varrho _{u\vert A^2}(W) + \varrho _{u\vert N}(W)] - [\delta _{l\vert A^1}(W) + \delta_{l\vert A^2}(W) + \delta _{l\vert N}(W)]
    \\ &= [\varrho _{u\vert A^2} (W) - \delta _{l\vert A^2}(W)] + [\varrho _{u\vert A^1}(W)-\delta _{l\vert A^1}(W)] + [\varrho _{u\vert N}(W) - \delta _{l\vert N}(W)]
    \\ &= b^*_2(X_2) - p^*_1(X_1) + \widetilde g(X_1-X_2) - \widetilde f(X_2-X_1),
  \end{split}
\end{align*}
which contradicts condition~\eqref{eq:gen1a}.

\medskip

{\bf Case 2.} \ $v_0^1, v_0^2 \notin W$.
Let $W^1 \coloneqq W\cap V^1$ and $W^2 \coloneqq W\cap V^2$, and let
\[
  X_1 \coloneqq \{ ({i,j}) \in S_{m,n} : v^1_{i,j} \in W^1 \} \ \ \ \hbox{and}\ \ \ X_2 \coloneqq \{ ( {i,j}) \in S_{m,n} : v^2_{i,j} \in W^2 \}.
\]
Similarly to Case 1, the contributions of $A^1$ and $A^2$ to $\varrho_u(W) - \delta_l(W)$ are
\begin{align*}
  \begin{split}
    &\varrho _{u\vert A^1}(W) - \delta _{l\vert A^1}(W) = \widetilde u(A^1 \cap \Delta^{\operatorname{in}}(W^1))  - \widetilde l(A^1 \cap \Delta^{\operatorname{out}} (W^1)) = b^*_1(X_1),\\
    &\varrho_{u\vert A^2}(W) - \delta _{l\vert A^2}(W) = \widetilde u(A^2 \cap \Delta^{\operatorname{in}}(W^2)) - \widetilde l(A^2 \cap \Delta^{\operatorname{out}} (W^2)) = - p^*_2(X_2) .
  \end{split}
\end{align*}
The arcs in $N$ contribute
\[
  \varrho _{u\vert N}(W) - \delta _{l\vert N}(W) = \widetilde g(X_2-X_1) - \widetilde f(X_1-X_2),
\]
and the special arc $a_0$ contributes $0$ as neither of its endpoints lies in $W$.

Summing all contributions, we get that
\begin{align*}
  \begin{split}
    0 &> \varrho _u(W) - \delta _l(W) = [\varrho _{u\vert A^1} (W) + \varrho_{u\vert A^2}(W) + \varrho _{u\vert N}(W)] - [\delta _{l\vert A^1}(W) + \delta _{l\vert A^2}(W) + \delta _{l\vert N}(W)]
    \\&= [\varrho _{u\vert A^1}(W) - \delta _{l\vert A^1}(W)] + [\varrho _{u\vert A^2}(W) - \delta_{l\vert A^2}(W) ] + [\varrho _{u\vert N}(W) - \delta _{l\vert N}(W)]
    \\&= b^*_1(X_1) - p^*_2(X_2) + \widetilde g(X_2-X_1) - \widetilde f(X_1-X_2),
  \end{split}
\end{align*}
contradicting~\eqref{eq:gen1b}.
\medskip

{\bf Case 3.} \ $v^1_0 \notin W, v^2_0 \in W$.
Let $W^1 \coloneqq V^
1\cap W$ and $W_2 \coloneqq V^2-W$, and let
\[
  X_1 \coloneqq \{ ( {i,j}) \in S_{m,n} : v^1_{i,j} \in W^1 \} \ \ \ \hbox{and}\ \ \ X_2 \coloneqq \{ ( {i,j}) \in S_{m,n} : v^2_{i,j} \in W^2 \}.
\]
The contribution of the arcs in $A^1$ and $A^2$ to $\varrho_u(W) - \delta_l(W)$ are
\begin{align*}
  \begin{split}
    &\varrho _{u\vert A^1}(W) - \delta _{l\vert A ^1}(W) = \widetilde u(A^1 \cap \Delta^{\operatorname{in}}(W^1)) - \widetilde l(A^1 \cap \Delta^{\operatorname{out}}(W^1)) = b^*_1(X_1),\\
    &\varrho _{u\vert A^2}(W) - \delta _{l\vert A^2}(W) = \widetilde u(A^2 \cap \Delta^{\operatorname{out}}(W^2)) - \widetilde l(A^2 \cap \Delta^{\operatorname{in}}(W^2)) = b^*_2(X_2).
  \end{split}
\end{align*}
The arcs in $N$ contribute
\[
  \varrho _{u\vert N}(W) - \delta _{l\vert N}(W) = \widetilde g(\ol{X_1} \cap \ol{X_2}) - \widetilde f(X_1\cap X_2).
\]
In this case, the special arc $a_0$ leaves $W$, thus its contribution to $\varrho_u(W)-\delta_l(W)$ is $-\alpha$.

Summing all contributions, we get that
\begin{align*}
  \begin{split}
    0 &> \varrho _u(W) - \delta _l(W) = -\alpha + [\varrho _{u\vert A^1}(W) - \delta _{l\vert A^1}(W)] + [\varrho _{u\vert A^2}(W) - \delta _{l\vert A^2}(W)]\\
      &\quad+ [\varrho _{u\vert N}(W) - \delta _{l\vert N}(W)] = -\alpha + b^*_1(X_1) + b^*_2(X_2) + \widetilde g(\ol{X_1}\cap \ol{X_2}) - \widetilde f(X_1\cap X_2),
  \end{split}
\end{align*}
contradicting~\eqref{eq:gen1alfa}.

\medskip

{\bf Case 4.} \ $ v^1_0 \in W, v^2_0 \notin W$.
Let $W^1 \coloneqq V^1-W$ and $W_2 \coloneqq V^2\cap W$, and let
\[
  X_1 \coloneqq \{ ( {i,j}) \in S_{m,n} : v^1_{i,j} \in W^1 \} \ \ \ \hbox{and}\ \ \ X_2 \coloneqq \{ ( {i,j}) \in S_{m,n} : v^2_{i,j} \in W^2 \}.
\]
The contribution of the arcs in $A^1$ and $A^2$ to $\varrho_u(W) - \delta_l(W)$ are
\begin{align*}
  \begin{split}
    &\varrho _{u\vert A^1}(W) - \delta _{l\vert A ^1}(W) = \widetilde u(A^1 \cap \Delta^{\operatorname{out}}(W^1)) - \widetilde l(A^1 \cap \Delta^{\operatorname{in}}(W^1)) = - p^*_1(X_1),\\
    &\varrho _{u\vert A^2}(W) - \delta _{l\vert A^2}(W) = \widetilde u(A^2 \cap \Delta^{\operatorname{in}}(W^2)) - \widetilde l(A^2 \cap \Delta^{\operatorname{out}}(W^2)) = -p^*_2(X_2).
  \end{split}
\end{align*}
The arcs in $N$ contribute
\[
  \varrho _{u\vert N}(W) - \delta _{l\vert N}(W) = \widetilde g(X_1\cap X_2) - \widetilde f(\ol{X_1} \cap \ol{X_2}).
\]
In this case, the special arc $a_0$ enters $W$, thus its contribution to $\varrho_u(W)-\delta_l(W)$ is $\beta$.

Summing all contributions, we get that
\begin{align*}
  \begin{split}
    0 &> \varrho _u(W) - \delta _l(W) = \beta + [\varrho _{u\vert A^1}(W) - \delta _{l\vert A^1}(W)] + [\varrho _{u\vert A^2}(W) - \delta _{l\vert A^2}(W)]\\
      &\quad+ [\varrho _{u\vert N}(W) - \delta _{l\vert N}(W)] = \beta - p^*_1(X_1) - p^*_2(X_2) + \widetilde g(X_1\cap X_2) - \widetilde f(\ol{X_1} \cap \ol{X_2}),
  \end{split}
\end{align*}
contradicting~\eqref{eq:gen1beta}, which completes the proof of Theorem~\ref{thm:general1}.
\FBOX

\medskip

We conclude by noting that a prefix-bounded matrix satisfying the prescribed entry bounds and total-sum bounds can be efficiently computed using standard network-flow algorithms, as detailed in Section~\ref{sec:arammodell}.
Furthermore, in the absence of such a matrix, the alternative proof of Theorem~\ref{thm:general1} --- based on Hoffman's circulation theorem --- provides a method for identifying the subsets $X_1,X_2 \subseteq S_{m,n}$ that certify infeasibility, in accordance with Theorem~\ref{thm:general1}.

\section*{Acknowledgments}

We are extremely grateful to Jack Edmonds who drew our attention to this topic and initiated an investigation to explore the link between ASMs and standard tools of combinatorial optimization.
We also deeply thank Richard Brualdi, Dylan Heuer, Karola M\'esz\'aros, Walter Morris, Cian O'Brien, Rachel Quinlan, and Jessica Striker for sharing their insights and perspectives with us.
We thank L\'aszl\'o V\'egh for his valuable comments on a previous version of the manuscript.

This research has been implemented with the support provided by the Ministry of Innovation and Technology of Hungary from the National Research, Development and Innovation Fund, financed under the ELTE TKP 2021-NKTA-62 funding scheme, and by the Ministry of Innovation and Technology NRDI Office within the framework of the Artificial Intelligence National Laboratory Program, and by the Ministry of Innovation and Technology of Hungary from the National Research, Development and Innovation Fund --- grant number ADVANCED~150556.

\bibliographystyle{plain_no_pages} 
\bibliography{bibliography}

\end{document}